\documentclass{article}
\usepackage{def}
\title{Reduction of the Hurwitz space of metacyclic covers}
\author{Irene I.\ Bouw}
\date{}
\begin{document}
\maketitle


\begin{abstract} We compute the stable reduction of some Galois covers
  of the projective line branched at three points. These covers are
  constructed using Hurwitz spaces parameterizing metacyclic covers.
  The reduction is  determined by a certain hypergeometric
  differential equation. This generalizes the result of
  Deligne--Rapoport on the reduction of the modular curve $X(p)$.
\end{abstract}

\section*{Introduction}
Recently, much progress has been made in understanding the reduction
of Galois covers of the projective line from characteristic zero to
characteristic $p>0$.  The starting point is a series of papers by
Raynaud, \cite{Raynaudfest}, \cite{Raynaud94}, \cite{Raynaud98} in
which the theory of the stable reduction of a Galois cover is
developed.  In \cite{Raynaud98}, one finds a criterion for good
reduction for Galois covers of the projective line branched at three
points, in case the order of the Galois group is strictly divisible by
$p$. Raynaud's Criterion implies, for example, that many covers which
have been constructed with the method of { rigidity} have good
reduction to characteristic $p$.

Raynaud's results have been extended in \cite{special} and
\cite{Stefan01}. Let $f:Y\to \PP^1$ be a $G$-Galois cover branched at
three points $x_1, x_2, x_3$, defined over a field $K$ of
characteristic zero.  Suppose that $p>2$ is a prime which strictly
divides the order of $G$. It is proved that the stable reduction of
$f$ at a prime of $K$ above $p$ is of a very easy kind.  This theorem
implies that to describe the stable reduction $\bar{f}$ of $f$, it
suffices to compute the associated deformation datum (Definition
\ref{defodatum}) for which there are only finitely many possibilities.
However, in general, it is not known how to compute $\bar{f}$.

In the literature, the stable reduction of a Galois cover with
bad reduction has been computed in very few  cases. Most
attention has been focused on cyclic covers of the projective line.
Coleman and McCallum \cite{CoMc} computed the reduction of cyclic
covers of the projective line branched at three points. In
\cite{Henrio}, \cite{Saidi00}, \cite{Claus}, \cite{Matignon01} the
reduction of $p$-cyclic covers of $\PP^1$ branched at $r\geq 3$ points
is studied. But even in this case the results are far from complete.

Other known results use the reduction of modular curves. For example,
let $X(N)$ be the modular curve parameterizing elliptic curves with
full level-$N$ structure. Then $f:X(2p)\to X(2)$ is a Galois cover
with Galois group $\PSL_2(p)$ branched at three points of order $p$.
Its stable reduction can be deduced from \cite{DelRap} and
\cite{KatzMazur}. 

In this paper, we compute the stable reduction of Galois covers of the
projective line branched at three points with Galois group
$\PSL_2(p)$, $\SL_2(p)$ or $\PGL_2(p)$. The construction of these
covers goes back to V\"olklein, \cite{Voelklein95}. Essentially, the
covers we consider arise as Hurwitz spaces $H({\bf a})$ parameterizing
covers of the projective line branched at four points.  More
precisely, choose $m|(p-1)$ and let $N$ be an extension of $\ZZ/m$ by
$\ZZ/p$. Then $H({\bf a})$ parameterizes $N$-Galois covers of $\PP^1$
branched at $0,\,1,\,\infty,\, \lambda$ with ramification of order
prime-to-$p$ described by the type ${\bf a}=(a_1, a_2, a_3, a_4)$,
Definition \ref{deftype}. Write $\pi({\bf a}):H({\bf
  a})\to\PP^1_\lambda$ for the cover which sends $g$ to the branch
point $\lambda$. The cover we consider is the Galois closure
$\varpi({\bf a})$ of $\pi({\bf a})$. It is expected that all covers of
the projective line which are branched at three points and whose
Galois group is $\PSL_2(p)$ or $\SL_2(p)$ arise via a construction
like this,  Section \ref{exasec}.

\bigskip\noindent{\bf Theorem} \ref{grthm}. {\sl The cover $\varpi({\bf
  a})$ has good reduction if and only if $a_1+a_2+a_3+a_4\neq 2m$.}

\bigskip In Section \ref{exasec} we illustrate that this theorem does
not follow from Raynaud's criterion for good reduction. For $m=2$ the
cover $\varpi({\bf a})$ is just the cover $X(2p)\to X(2)$.

In case the cover $\varpi({\bf a})$ has bad reduction, we compute its
stable reduction in Section \ref{defodatumsec}. We do not give a
modular interpretation for the stable model, therefore our approach
differs from the approach of \cite{DelRap} and \cite{KatzMazur}.
However, also in our approach, the key ingredient for computing the
reduction of $\varpi({\bf a})$ is the study of the reduction of the
metacyclic covers which the Hurwitz space parameterizes. It turns out
that the stable reduction is determined by a so-called ${\bf a}$-Hasse
invariant $\Phi_{\bf a}$. For $m=2$ this is just the classical Hasse
invariant. An important property of $\Phi_{\bf a}$ is that it is the
solution to a hypergeometric differential equation, Section
\ref{Hassesec}. This hypergeometric differential equation determines
all combinatorial invariants which are associated to the stable
reduction, see Section \ref{defodatumsec} for a precise statement.

The organization of this paper is as follows. In Section
\ref{stablesec} we define the stable reduction and collect some
results on its structure. In Section \ref{metasec} we introduce the
metacyclic covers and study their reduction. In Section \ref{Hurwitz0}
we define the Hurwitz space $H({\bf a})$ in characteristic zero and
describe the map $\pi({\bf a}):H({\bf a})\to\PP^1_\lambda$ and its
Galois closure. In Section \ref{grsec} we prove the good reduction
theorem. Sections \ref{Hassesec} and \ref{defodatumsec} study the
stable reduction of the covers with bad reduction. First we define the
${\bf a}$-Hasse invariant, Section \ref{Hassesec}. Then we show how
this gives the stable reduction, Section \ref{defodatumsec}. The rest
of the paper contains complements and examples.

\section{The stable model}\label{stablesec}
In this section we recall the definition of the stable model. We refer
to \cite{Stefan01} for proofs and more details. Suppose that $G$ is a
finite group whose order is strictly divisible by $p$.  Let $R$ be a
complete discrete valuation ring with fraction field $K$ of
characteristic zero and residue field $k=\bar{k}$ of characteristic
$p>2$.  Let $f:Y\to X=\PP^1_K$ be a $G$-Galois cover branched at ${\bf
  x}=(x_1, \ldots, x_r)$, where the $x_i$ are distinct $K$-rational
points.  Write $\BB_0=\{1, \ldots, r\}$ for the set indexing the
branch points.  We assume that $r\geq 3$. We assume moreover that $(X;
{\bf x})$ has {\sl good reduction}, i.e.\ there exists a model
$X_{0,R}=\PP^1_R$ of $X$ over $R$ such that the $x_i$ extend to
pairwise disjoint sections $\Spec(R)\to X_{0,R}$.

Denote the ramification points of $f$ by ${\bf y}=(y_1, \ldots y_s)$.
Let $(Y_R; {\bf y})$ be the unique extension of $(Y; {\bf y})$ to a
stably pointed curve, which exists after replacing $R$ by a finite
extension, if necessary. Recall that $(Y_R;{\bf y})$ is called {\sl
  stably pointed} if $Y_R$ is a semistable curve such that the usual
``three point condition'' is fulfilled for the union of the marked
points and the singular points of $\Yb:=Y_R\otimes_R k$. The action of
$G$ extends to $Y_R$; denote the quotient of $Y_R$ by $G$ by $X_R$.
Then $f_R:Y_R\to X_R$ is called the {\sl stable model} of $f$.  The
special fiber $\bar{f}:\Yb\to \Xb$ of $f_R$ we call the {\sl stable
  reduction} of $f$.  We say that $f$ has {\sl good reduction} if and
only if the stable reduction $\bar{f}:\Yb\to \Xb$ is a separable
morphism. This is equivalent to $X_R$ being smooth. If $f$ does not
have good reduction, we say it has {\sl bad reduction}.

Suppose that $f$ has bad reduction.  Let $\fb:\Yb\to \Xb$ be the
stable reduction of $f$. We denote the dual graph of $\Xb$ by $T$. The
vertices $\VV$ of $T$ correspond to the irreducible components of
$\Xb$ and the (oriented) edges $\EE$ to the singular points. The
source and target of an edge $e$ are denoted by $s(e)$ and $t(e)$.

Since $(X; {\bf x})$ has good reduction, there exists a natural map
$X_R\to X_{0, R}$. This map restricts to the identity on a unique
irreducible component of $\Xb$ and contracts all other components of
the special fiber to a point.  Therefore we may (and will) consider
$\Xb_0$ as an irreducible component of $\Xb$, which we call the {\sl
  original component.} The dual graph $T$ of $\Xb$ is a tree.  We
denote the vertex of $T$ corresponding to the original component
$\Xb_0$ by $v_0$ and consider $T$ to be oriented from $v_0$.

An edge $e$ is called {\sl positive} if its orientation coincides with
the orientation of $T$.  For every $v\in \VV$, we denote the
corresponding component of $\Xb$ by $X_v$.  We choose a component
$Y_v$ of $\Yb$ above $X_v$. For every $v\in \VV$ we denote by
$\infty_v$ the unique point of $X_v$ which is singular in $\Xb$ and
such that the corresponding edge $e$ with $s(e)=v$ is negative. If $e$
is an edge with $s(e)=v$, we write $x_e$ for the singular point of
$X_v$ corresponding to $e$.

An irreducible component $X_v$ of $\Xb$ is called a {\sl tail} if it
is different from the original component and meets the rest of $\Xb$
in just one point.  We denote by $\BB$ the subset of $\VV$
corresponding to the tails of $\Xb$.  A tail is called {\sl primitive}
if one of the branch points specializes to this component; otherwise
the tail is called {\sl new}.  Let $\Bp$ (resp.\ $\Bn$) be the subset
of $\BB$ corresponding to the primitive (resp.\ new) tails. We
identify $\Bp$ with a subset of $\BB_0=\{1, \ldots, r\}$: we write
$X_i$ for the primitive tail to which the branch point $x_i$
specializes. The complement $\Bw:=\BB_0-\Bp$ corresponds to the branch
points whose ramification index is divisible by $p$. We use the
notation $x_i$ to denote both the point on $X$ and its specialization
to $\Xb$.

For every $b\in \BB$, the map $Y_b\to X_b$ is separable.  The inertia
group of a point above $\infty_b$ has order $pm_b$, for some $m_b$
prime-to-$p$.  Denote the corresponding conductor by $h_b$. We define
the {\sl ramification invariant} to be $\sigma_b={h_b}/{m_b}.$ (This
is the jump in the filtration of the higher ramification groups in the
upper numbering.)  If $b\in \Bn$, the cover $Y_b\to X_b$ is unbranched
outside $\infty_b$.  If $i\in \Bp$, then $Y_i\to X_i$ is tamely
branched at $x_i$ and unbranched outside $x_i$ and $\infty_i$.  If
$v\in\VV- \BB$, then the restriction of $\fb$ to $X_v$ is inseparable.

Choose an irreducible component $\Yb_0$ of $\Yb$ above the original
component $\Xb_0$. Since $f$ has bad reduction, the inertia group
$I_0$ if $\Yb_0$ has order $p$. The decomposition group $D_0$ of
$\Yb_0$ is the extension of a group $H_0$ of order prime-to-$p$ by
$I$. The map $\Yb_0\to\Xb_0$ is the composition of a purely
inseparable map $\Yb_0\to\Zb_0$ and an $H_0$-Galois cover
$\Zb_0\to\Xb_0$. Write $\Yb_0'$ (resp.\ $\Zb_0'$) for the quotient of
$\Yb_0$ (resp.\ $\Zb_0$) by the prime-to-$p$ part of the center of
$D_0$. Then $\Zb_0'\to\Xb_0$ is an $m$-cyclic cover, for some $m$
dividing $p-1$. Let $\chi:\ZZ/m\to\FF_p^\times$ be an injective
character. We choose generators $\phi$ and $\psi$ of the decomposition
group $D_0'$ of $\Yb_0'$ such that $\phi^p=\psi^m=1$ and
$\psi^i\phi\psi^{-i}=\phi^{\chi(i)}.$

The restriction of $\Yb_0\to \Zb_0$ to some open $U$ of $\Zb_0$ is an
$\alphap$ or $\mup$-torsor. We say that $f$ has {\sl multiplicative
  (bad) reduction} if $\Yb_{0}\to \Zb_{0}$ restricts to a
$\mup$-torsor. If $\Yb_{0}\to \Zb_{0}$ restricts to an
$\alphap$-torsor, we say that $f$ has {\sl additive (bad) reduction.}

To an $\alphap$ or $\mup$-torsor, we can associate a differential
$\omega$, \cite[Proposition 4.14]{Milne}.  Let $W$ be the inverse
image of $U$ in $\Yb_0$.  There exists an open covering $\{U_j\}_j$ of
$U$ and elements $f_j\in\Gamma(U_j, \O)$ such that the torsor
restricted to $U_j$ is given by an equation $v^p=f_j$.
 
 Suppose first that $W\to U$ is a $\mup$-torsor, then $f_j\in
 \Gamma(U_j, {\cal O})^\ast$.  To the torsor we associate the
 differential
\[\omega=\frac{{\rm d} f_j}{f_j}\in H^0(U,\Omega^1) .\]
Here we assume that we fixed an isomorphism $\ZZ/p\simeq \mup(K)$.
Note that $\omega$ does not depend on $j$ or the choice of $f_j$.  In
case $W\to U$ is an $\alphap$-torsor, we associate to it the
differential
\[\omega={{\rm d} f_j}\in H^0(U,\Omega^1).\]
Moreover, $\psi^\ast\omega=\chi(1)\omega$. See \cite{Stefan01} for
more details.

\begin{defn}\label{defodatum} We call the tuple $(\Zb_0', \omega)$ the
  {\sl deformation datum} corresponding to $f$.
\end{defn}

\begin{lem}\label{polesomega}
\begin{itemize}
\item[(a)] The differential $\omega$ is regular outside $x_b$ for
  $b\in \Bw$. For $b\in \Bw$, we have $\Ord_{z_b}(\omega)=-1$. Here
  $z_b\in\Zb_0$ is a point above $x_b$.
\item[(b)] The differential $\omega$ has a zero at $z\in \Zb_0$ if and
  only if the image of $z$ in $\Xb_0$ is a singular point of $\Xb$.
\end{itemize}
\end{lem}

\proof First suppose that $w\in\Zb_0$ is a point whose image in
$\Xb_0$ is $x_b$ for $b\in\Bw$. Let $t$ be a local parameter at $w$.
Then $\Yb_0\to \Zb_0$ is given by $y^p=t^n \bmod t^{n+1},$ locally
around $w$, for some $n>0$ which is prime-to-$p$. It follows that
$\omega$ is a logarithmic differential which has a simple pole at $w$.

Let $w\in \Zb_0$ be a point which does not map to $x_b$ for $b\in\Bw$.
Let $t$ be a local parameter at $w$. Then the natural map
$\Yb_0\to\Zb_0$ is given by $y^p=1+t^n \bmod t^{n+1},$ locally around
$w$, \cite[p.\ 192]{Raynaudfest}.  It follows from the definition of
$\omega$ that $\omega$ does not have a pole at $w$ and that $\omega$
has a zero at $w$ if and only if $w$ is singular in $\Zb$. \Endproof

\begin{defn} \label{conddef}
  Let $\EE_0=\{e\in \EE \, |\, s(e)=v_0\}$. Let $z_e\in\Zb_0$ be a point
  above $x_e$. For $e\in\EE_0$, the {\sl conductor} of $e$ is defined
  as
\[h_e=\Ord_{z_e}(\omega)+1. \]
\end{defn}

Suppose that $e\in\EE_0$. Let $m_e$ be the
ramification index of $z_e$ in $\Zb_0\to\Xb_0$. Define
\[\nu_e:=\left[\frac{h_e}{m_e}\right], \quad a_e:=m\left(\frac{h_e}{m_e}-\nu_e\right).\]

\begin{lem}\label{vc} Let $e\in\EE_0$.
\begin{itemize}
\item[(a)] If the subtree of $T$ with root $t(e)$ contains a new tail,
  then $\nu_e\geq 1$.
\item[(b)] 
\[\sum_{e\in\EE_0} \nu_e=|\EE_0|-2-\frac{1}{m}\sum_{b\in\BB} a_b.\]
\end{itemize}
\end{lem}

\proof The lemma is proved in \cite{Stefan01}. Part (a) is a
reformulation of \cite[Proposition 3.3.5]{Raynaud98}. Part (b) follows
from the vanishing cycle formula \cite[Section 3.4.2]{Raynaud98}.
\Endproof

\label{3ptsec}

\begin{thm}\label{specialthm}
  Let $G$ be a group whose order is strictly divisible by $p$.
  Let $f:Y\to X=\PP^1_K$ be a $G$-Galois cover over $K$
  branched only at $0,1,\infty$. Suppose that $f$ has bad
  reduction. 
\begin{itemize}
\item[(a)] The cover  $f$ has multiplicative reduction.
\item[(b)] Every irreducible component of $\Xb$ other than $\Xb_0$ is a
  tail.
\item[(c)] Let $e\in\EE_0$ with $t(e)=b\in\BB$. Choose a point $z_e$
  of $\Zb_0'$ above $x_e$. Then $z_e$ is a ramification point of
  $\Zb_0'\to\Xb_0$. Write $m_e$ for its ramification index.  The
  ramification invariant $\sigma_b$ of the tail $X_b$ is equal to
  $h_e/m_e$.
\end{itemize}\end{thm}

\proof \cite{Stefan01}.
\Endproof

\begin{rem}\label{typehurw} We use the notation of Theorem
  \ref{specialthm}.  Lemma \ref{vc} and Theorem \ref{specialthm} imply
  that $\nu_e=1$ if $t(e)\in\Bn$ and $\nu_e=0$ if $t(e)\in\Bp$.
  Moreover, we have that $\sum_{b\in\BB} a_b=m$. One easily checks
  from the proof of Lemma \ref{polesomega} that $\psi^{a_b}$ is the
  canonical generator of the inertia group of $z_e$ in
  $\Zb_0'\to\Xb_0$ (with respect to the character $\chi$). In the
  terminology of Definition \ref{deftype}, we say that the {\sl type}
  of $\Zb_0'\to\Xb_0$ is $(x_e, a_b)$, where $e\in\EE_0$ and $b=t(e)$.
\end{rem}

\section{The reduction of metacyclic covers}\label{metasec}
In this section, we study metacyclic covers of the projective line
branched at four points. We start by defining the class of metacyclic
covers we are interested in in this paper. We determine the different
types of reduction that can occur. In Definition \ref{casedef}, we
distinguish the three different cases that may occur. These play an
important role in the rest of the paper.

\label{setup}
We fix the following data.
\begin{itemize}
\item Let $p\neq 2$ be a
 prime number and $m>1$ an integer such that $p\equiv 1\bmod{m}$.
\item Let $R$ be a complete discrete valuation ring whose fraction
  field $K$ has characteristic zero and whose residue field $k$ is
  algebraically closed of characteristic $p$.
\item Fix a character $\chi:\ZZ/m\to \FF_p^\times$ of order $m$.  Put
  $\zeta_m=\chi(1)$. We choose a lift $\chi_K:\ZZ/m\to K^\times$ of
  $\chi$. We denote $\chi_K$ also by $\chi$.
\item Let $N$ be an extension of $\ZZ/m$ by $\ZZ/p$.  We fix once and
  for all generators $\phi$ and $\psi$ of $N$. We suppose that
  $\phi^p=\psi^m=1$ and $\psi^i\phi\psi^{-i}=\phi^{\chi(i)}$.
\end{itemize}

\begin{defn}\label{deftype}
Let ${\bf x}=(x_1, \ldots x_r)$ be $r$ distinct points of
$X=\PP^1_{\bar{K}}$. Let ${\bf a}=(a_1, \ldots, a_r)$ be an $r$-tuple of integers with $0< a_i<m$
and $\sum a_i\equiv 0\bmod{m}.$ Let $f:Y\to X$ be an $N$-Galois cover.
We say that $f$ is {\sl of type }${\bf (x; a)}$ (with respect to
$\chi$) if the following holds.
\begin{itemize}
\item[(a)] The cover $f$ is  branched only at $x_1, \ldots, x_r$.
\item[(b)] The element $\psi^{a_i}$ is the canonical generator with
  respect to $\chi$  of some
  point $y_i$ of $Y$ above $x_i$. In other words, if $u_i$ is a local
  parameter of $y_i$ then
\[  (\psi^{a_i})^\ast u_i\equiv \chi(1)^{(a_i, m)} u_i \pmod{u_i^2}.\]
\end{itemize}
 
\end{defn}

In concrete terms this means the following. Let $f:Y\to X$ be a
metacyclic cover of type ${\bf (x, a)}$ and suppose that none of the
$x_i$ is equal to $\infty$, for simplicity. Let $Z$ be the quotient of
$Y$ by the normal subgroup of $N$ of order $p$. Then $Y\to Z$ is
\'etale and $Z$ is the complete nonsingular curve given by the
equation
\[z^m=\prod_i (x-x_i)^{a_i},\]
where $\psi^\ast z=\chi(1) z.$
We  also use the terminology `of type ${\bf( x; a)}$' for the (unique)
$m$-cyclic cover $Z\to X$.

\begin{defn}\label{xidef}
  Let $g:Z\to X$ be the $m$-cyclic cover of type ${(\bf x; a)}$. Write
  $V:=\Hom(\pi_1(Z), \ZZ/p)$. Then $\ZZ/m$ acts on $V$ and we may
  consider $ V_\chi:=\{\xi\in V\,|\, \psi\xi=\chi(1)\xi\}.$ Suppose
  $f:Y\to X$ is a metacyclic cover of type ${(\bf x; a)}$. Then $f$
  factors through $g$. We define $\xi_f\in V$ as the element
  corresponding to the exact sequence
\[ 1\to \pi_1(Y)\to \pi_1(Z)\to\Gal(Y, Z)=\ZZ/p\to 1.\]
One checks that $\psi\xi_f=\chi(1)\xi_f$ and therefore $\xi_f\in V_\chi$.
\end{defn}

\begin{lem}\label{dimchar0}
With notations as above:
\[\dim_{\FF_p} V_\chi=r-2.\]
\end{lem}

\proof This follows from the main result of \cite{CW}. See also
\cite[Section 2.4]{Voelklein95}.  \Endproof

Two metacyclic covers $f_i:Y_i\to X$ of type ${\bf( x; a)}$ are
isomorphic   if there exists an $N$-equivariant
isomorphism $h:Y_1\to Y_2$ such that $f_1=f_2\circ h$. Note that such
$h$, if it exists, is unique, since the center of $N$ is trivial. Two
metacyclic covers $f_1, f_2$ of type ${\bf (x; a)}$ are isomorphic if
and only if the corresponding elements $\xi_i\in V_\chi$ are in the
same orbit under the action of $\ZZ/m$ on $V_\chi$. Therefore Lemma
\ref{dimchar0} implies that, up to isomorphism, there are
$(p^{r-2}-1)/m$ metacyclic covers of $X$ of type ${\bf (x;
  a)}$.

Since we suppose that $(X; {\bf x})$ has {\sl good reduction}, the
$m$-cyclic cover $Z\to X$ extends to an $m$-cyclic cover $Z_{0,R}\to
X_{0,R}$ of smooth curves branched  at the $x_i$.  Let $\Zb_0$ and
$\Xb_0$ be the special fibers of $Z_{0, R}$ and $X_{0, R}$.   Let
$\bar{V}=\Hom(\pi_1(\Zb_0), \ZZ/p)$ and let $\bar{V}_{\chi}$ be the
$\chi$-isotypical part of $\bar{V}$. Since every \'etale $p$-cyclic
cover of $\Zb_0$ can be uniquely lifted to an \'etale $p$-cyclic cover
of $Z$, we obtain a canonical injection $\bar{V}_\chi\hookrightarrow
V_\chi$. This proves the following criterion for good reduction.

\begin{lem}\label{goodredlem}\label{goodredequal}
  Let $f:Y\to X$ be a metacyclic cover of type ${{\bf (x; a)}}$ and
  let $\xi\in V_\chi$ be the corresponding element. Then $f$ has
  potentially good reduction if and only if $\xi\in
  \bar{V}_\chi$.
\end{lem}

\label{altformsec} 

Artin--Schreier Theory implies that there exists a canonical
isomorphism $\bar{V}\simeq H^1(\Zb, {\cal O}_{\Zb})^F$, where $F$
denotes the absolute Frobenius.  Write
\[ H^1(\Zb_0, {\cal O}) =\oplus_{i=1}^{m-1} L_{\chi^i},\]
where $L_{\chi^i}=\{\xi|\psi\xi=\chi^i(1)\xi\}$. Then
$\bar{V}_\chi=L_\chi^F$.

\begin{lem}\label{dimLi}
The dimension of $L_\chi$ as a $k$-vector space is $(a_1+\cdots
+a_r-m)/m$.
\end{lem}

\proof
\cite[Lemma 4.3]{thesis}.
\Endproof

Note that $0\leq \dim_k L_\chi\leq r-2$. In fact, $\dim_k
L_\chi+\dim_k L_{\chi^{-1}}=r-2$.  This gives an upper bound on the
number of metacyclic covers of type ${\bf (x; a)}$ in characteristic
$p$. If $\dim L_\chi <r-2$ there are strictly less $N$-covers of type
${(\bf x;a)}$ in characteristic $p$ than in characteristic zero. In
principle one can compute $\dim_{\FF_p} \bar{V}_\chi$ by using the
explicit basis for $L_\chi$ and the matrix of $F$ with respect to this
basis as computed in \cite[Section 5]{thesis}.

\label{duality}
Suppose that $f:Y\to X$ is a metacyclic cover of type ${\bf (x;a)}$
with bad reduction.  Lemma \ref{polesomega} implies that $\omega$ is a
regular differential on $\Zb_0$, i.e.\ an element of $H^0(\Zb_0,
\Omega^1)$. The vector space $H^0(\Zb_0, \Omega^1)$ is dual to
$H^1(\Zb_0, {\cal O}_{\Zb_0})$. The transpose of the absolute
Frobenius $F:H^1(\Zb_0, {\cal O}_{\Zb_0})\to H^1(\Zb_0, {\cal
  O}_{\Zb_0})$ is called the Cartier operator $\C: H^0(\Zb_0,
\Omega^1)\to H^0(\Zb_0, \Omega^1)$, \cite{SerreMex}. An element
$\omega\in H^0(\Zb_0, \Omega^1)$ defines a $\mup$-torsor if
$\C\omega=\omega$ and an $\alphap$-torsor if $\C\omega=0$,
\cite[Proposition 4.14]{Milne}. Let $\bar{V}=H^1(\Zb_0, {\cal
  O}_{\Zb_0})$ and $\bar{V}_\chi=\{\xi\in \bar{V}|\,
\psi\xi=\chi(1)\xi\}$. Note that this definition coincides with the
definition we gave in Section \ref{setup}.

Define $\bar{M}=H^0(\Zb_0, \Omega^1)^{\C}$ and let
$\bar{M}_\chi=\{\omega\in \bar{M}\,|\, \psi\omega=\chi(1) \omega\}$.
Note that $\bar{M}_\chi$ is an $\FF_p[\ZZ/m]$-module.  Since $\C:
H^0(\Zb_0, \Omega^1)\to H^0(\Zb_0, \Omega^1)$ is the transpose of
$F:H^1(\Zb_0, {\cal O}_{\Zb_0})\to H^1(\Zb_0, {\cal O}_{\Zb_0})$, we
conclude that
\[\bar{M}_\chi=\bar{V}_{\chi^{-1}}^{\dual},\]
where `dual' means duality of $\FF_p[\ZZ/m]$-modules. 

\begin{prop}\label{M1quotient}
There exists a surjective morphism of $\FF_p[\ZZ/m]$-modules
\begin{equation}V_\chi\twoheadrightarrow \bar{M}_\chi.\label{M1quoeq}\end{equation}
\end{prop}

\proof Since $\bar{M}_\chi=\bar{V}_{\chi^{-1}}^{\dual}$, there is an
inclusion $\bar{M}_{\chi}^{\dual}\hookrightarrow V_{\chi^{-1}}$. 
The Weil pairing defines a map
\[V_\chi\times V_{\chi^{-1}}\to\mup.\]
The choice of an isomorphism $\ZZ/p\simeq \mup(K)$ identifies
therefore $V_\chi$ and $V_{\chi^{-1}}^{\dual}$. This defines
(\ref{M1quoeq}).\Endproof

In more concrete terms we can describe the surjection of Proposition
\ref{M1quotient} as follows. Let $f:Y\to X$ be a metacyclic cover of
type ${{\bf (x; a)}}$ and let $\xi$ be the corresponding element of
$V_\chi$. If $f$ has multiplicative reduction, then $\xi$ maps to
$\omega_0$ in $\bar{M}_\chi$. Otherwise, $\xi$ maps to zero. In case
$H^0(\Zb, \Omega_1)^{\C=0}=0$ there are no metacyclic covers with
additive reduction and the following sequence is exact:
\[0\to \bar{V}_\chi\to V_\chi\to \bar{M}_\chi\to 0.\]
This happens, for example, when $\Zb_0$ is ordinary. Namely, in that
case $F:L_{\chi^i}\to L_{\chi^i}$ is a bijection for every $i$, in
particular for $i=1, -1$. Therefore $\dim_{\FF_p} \bar{V}_\chi=\dim_k
L_\chi$ and $\dim_{\FF_p} \bar{M}_\chi=\dim_k L_{\chi^{-1}}$.
Therefore $H^0(\Zb, \Omega_\chi)^{\C=0}=0$.

\label{4pointssec}
Suppose that $r=4$.  Define $a:=\dim_{\FF_p} \bar{V}_\chi=\dim_{\FF_p}
L_\chi^F$ and
$b:=\dim_{\FF_p}\bar{M}_\chi=\dim_{\FF_p}L_{\chi^{-1}}^F$. Note that
$a,b\leq 2$. Lemma \ref{dimLi} implies that $\dim_k
L_{\chi^i}=r-2=2$, therefore $a,b\leq 2$.  The following proposition
lists the possibilities for $(a,b)$.

\begin{prop}\label{pdivisible}
  Suppose $r=4$. Then $(a,b)\in\{(2,0), (0,2), (1,1), (0,0)\}$.
\end{prop}

\proof The $\FF_p[\ZZ/m]$-module $\bar{V}$ is equal to $\Hom (\mup,
J(\Zb_0)[p])$, where `Hom' should be considered in the category of
finite flat group schemes.  The $p$-divisible group
$J(\Zb_0)[p^\infty]$ decomposes into eigenspaces of the automorphism
$\psi$, since $\ZZ_p$ contains the $m$th roots of unity. Write
\[J(\Zb_0)[p^\infty]=\oplus_i J(\Zb_0)[p^\infty]_{\chi^i}\]
for this decomposition.  We conclude that the finite flat group scheme
$J(\Zb_0)[p]_{\chi}= J(\Zb_0)[p^\infty]_{\chi}/p$ is a direct summand of
the group scheme $J(\Zb_0)[p]$. Recall that there is a direct sum
decomposition
\[J(\Zb_0)[p]= (\ZZ/p)^\sigma \times (\mup)^\sigma\times \Lambda,\]
where $\sigma$ is the $p$-rank of $\Zb_0$ and $\Lambda$ is some
local-local finite flat group scheme. Since $J(\Zb_0)$ is an abelian
variety, $J(\Zb)[p]$ does not have any direct summand $\alphap$.

Since the dimension of $J(\Zb_0)[p]_{\chi}$ is two, we conclude that
it is isomorphic to either $(\ZZ/p)^2$, $(\mup)^2$, $\ZZ/p\oplus \mup$
or it is a local-local group scheme.  The definition of $a$ and $b$
implies that
\[ J(\Zb_0)[p]_{\chi}\simeq (\ZZ/p)^{b}\times (\mup)^a\times \Lambda_{\chi},\]
where $\Lambda_{\chi}$ is local-local.  This proves the proposition.
\Endproof

It is easy to check that all possibilities listed in Proposition
\ref{pdivisible} occur; we give explicit conditions for each
possibility to occur further on.  

\begin{defn}\label{casedef}
Distinguish three cases.
\begin{itemize}
\item {\bf  The multiplicative case}:  $a_1+a_2+a_3+a_4=m$.
\item {\bf  the mixed case}:   $a_1+a_2+a_3+a_4=2m$.
\item {\bf the \'etale case}: $a_1+a_2+a_3+a_4=3m$.
\end{itemize}
\end{defn}

\begin{prop} \label{mindim}
  Suppose $a_1+a_2+a_3+a_4=m$ and $\lambda\not\equiv
  0,1,\infty\bmod{p}$. Then all metacyclic covers of type $({\bf
    x;a})$ have multiplicative bad reduction. Therefore $(a,b)=(0,2)$, for
  all $\lambda\in\PP^1_k-\{0,1,\infty\}$.
\end{prop}

\proof The proposition is proved in \cite[Proposition 1.3]{special}
for every $r\geq 3$. \Endproof

\begin{prop} \label{genord}
Suppose $a_1+a_2+a_3+a_4=2m$. Then  $F:L_\chi\to L_\chi$ is an
isomorphism, for sufficiently general $\bar{\lambda}\in k$.
\end{prop}

\proof
  In fact, a similar statement holds without assumption on $r$ and the
  type. In \cite[lemma 4.6]{thesis} it is shown that the locus $U$ where
  $F:L_\chi\to L_\chi$ is an isomorphism is open. The fact that $U$ is
  non-empty follows from \cite[Proposition 7.4]{thesis}.  \Endproof

  Let $U\subset\PP^1_k-\{0,1,\infty\}$ be the locus where $F:L_\chi\to
  L_\chi$ is an isomorphism. In case $m=2$, this $U$ is the subset of
  $\lambda\in \PP^1_k-\{0,1,\infty\}$ such that the elliptic curve
  with Weierstrass equation $y^2=x(x-1)(x-\lambda)$ is ordinary.
  Therefore, in the mixed case, $U$ is in general strictly smaller than
  $\PP^1_k-\{0,1,\infty\}$, in contrast to what happens in the
  \'etale case (see below).
  
  For $\lambda\in U$ we have $(a,b)=(1,1)$.  Lemma \ref{goodredlem}
  implies that $(p-1)/m$ metacyclic cover of type $({\bf x;a})$ have
  good reduction and the other $p(p-1)/m$ have bad reduction. There
  are finitely many ${\lambda}\in \PP^1_k$ for which $F:L_\chi\to
  L_\chi$ is the zero morphism. For these ${\lambda}$ we have
  $(a,b)=(0,0)$ and all metacyclic covers of type $({\bf x;a})$ have
  additive bad reduction.

\begin{prop}\label{maxdim}  Let $f:Y\to X$ be a metacyclic cover of
  type ${{\bf (x; a)}}$.  Suppose that $a_1+\cdots+a_4=3m$ and
  $\lambda\not\equiv 0,1,\infty\bmod{p}$. Then $f$ has good reduction
  and $(a,b)=(2,0)$.
\end{prop}

\proof Let $f$ as in the statement of the proposition.  The assumption
$a_1+\cdots+a_4=3m$ implies that $\dim_k L_\chi=2$. Since $\dim_k
L_\chi+\dim_k L_{\chi^{-1}}=r-2$, we conclude that $\dim_k
L_{\chi^{-1}}=0$. Therefore Proposition \ref{mindim} implies that
$V_{\chi^{-1}}=\bar{M}_{\chi^{-1}}$.  By duality, we conclude that
$V_{\chi}=\bar{V}_{\chi}$. Therefore $f$ has good reduction.
\Endproof

\section{Description of the Hurwitz space in characteristic zero}\label{Hurwitz0}

In this section, we define the Hurwitz space parameterizing metacyclic
covers in characteristic zero. 

Let $G$ be a group whose order is strictly divisible by $p$.  Let
$R_0=W(\bar{\FF}_p)$ and $K_0$ its quotient field.  Let $N\simeq
\ZZ/p\rtimes \ZZ/m$ be as in Section \ref{setup}.  Recall that we have
fixed generators $\phi$ and $\psi$ of $N$ which satisfy
\[\phi^p=\psi^m=1, \qquad \psi\phi\psi^{-1}=\phi^{\chi(1)},\]
where $\chi:\ZZ/m\to\FF_p^\times$ is an injective character.  We
represent the element $\phi^i\psi^j$ of $N$ as $[i, \chi(1)^j]$, with
$i\in \FF_p$ and $\chi(1)^j\in \Im(\chi)\subset \FF_p^\times$.

Fix a type ${{\bf a}}$ with $\gcd(m, a_1, \ldots, a_4)=1$. The inverse
image of $b\in \ZZ/m$ is a conjugacy class of $N$ which we denote by
$C(b)$.  Let ${\bf C}=(C({a_1}), C({a_2}), C({a_3}), C({a_4}))$. Put
$\zeta_i:=\chi(1)^{a_i}\in \FF_p^\times$. Let $\QQ_p({\bf a})$ be the
smallest extension of $\QQ_p$ over which ${\bf C}$ is rational,
\cite[Section 7.1]{Serretopics}.

\begin{defn}
  Define $H({\bf a})=\hin_{4}({\bf a})/ \QQ_p({\bf a})$ as the (inner)
  Hurwitz space parameterizing $N$-Galois covers $Y\to X=\PP^1$ of
  type ${\bf (x;a)}$.  Write $\pi({\bf a}):H({\bf a})\to\PP^1_\lambda
  -\{0,1,\infty\}$ for the cover defined by $[f]\mapsto \lambda$.  We
  denote th Galois closure of $\pi({\bf a})$ by $\varpi({\bf a
    }):\HH({\bf a})\to\PP^1_\lambda$.
\end{defn}

If ${\bf a}$ is understood, we sometimes drop the index ${\bf a}$ from
the notation.  For generalities on Hurwitz spaces, we refer to
\cite{diss}.  This Hurwitz space has also been considered in
\cite{Berger}. Berger considers types of the form ${\bf
  a}=(1,1,1,m-3)$. He calls the Hurwitz space $H({\bf a})$ a ``fake
modular curve'' and shows that it is a quotient of the complex upper
half-plane by a non-congruence subgroup.

\begin{lem}
 The degree of $\pi\otimes {\bar{\QQ}_p}$ is $(p^2-1)/m$.
\end{lem}

\proof
 Let $f:Y\to X$ be an $N$-cover of type ${{\bf (x; a)}}$ and $Z\to X$ the
intermediate $m$-cyclic cover. The \'etale $p$-cyclic cover $Y\to Z$
corresponds to a nonzero element $\xi$ in $V_\chi\simeq (\ZZ/p)^2$,  Lemma
\ref{dimchar0}.  Two elements $\xi_1$ and $\xi_2$ define the same
cover if $\psi^i(\xi_1)=\xi_2$, for some $i$.  \Endproof

 One can describe describe $\hb\otimes \bar{\QQ}_p$ using Nielsen
 tuples, see for example \cite[Section 10.1.7]{Voelklein}.  Define
\[ \E({\bf a})=\{{\bf g}=(g_1, g_2, g_3, g_4)| \,g_i\in C(a_i),\,  N=\langle g_i \rangle,\,  \prod_i g_i=1\, \}.\]
The  Nielsen class is equal to
\[\Niin_4({\bf a})=\E({\bf a})/N,\]
where the group $N$ acts via uniform conjugation. The elements of this
set are called the {\sl Nielsen tuples}. Choose a presentation of the
fundamental group
\[\pi_1(\PP^1-{\bf x}, \ast)=\langle\gamma_1, \ldots, \gamma_4|\prod_i\gamma_i=1\rangle,\]
where $\gamma_i$ corresponds to a ``loop'' around the point $x_i$, as
usual. One can identify the  Nielsen tuples with
isomorphism classes of surjective homomorphisms
\[\pi_1(\PP^1-{\bf x}, \ast)\to N, \qquad \Im(\gamma_i)\in C(a_i). \]
We also consider $\E^\ast({\bf a}):=\E({\bf a})/\FF_p$, where
$\FF_p\subset N$ acts via uniform conjugacy.

Choose a base point $\lambda_0\in\PP^1-\{0,1,\infty\}$.  The cover
$\pi:H\to \PP^1_{\lambda}$ defines an action of $\Pi:=\pi_1(\PP^1-\{0,
1, \infty\}, \lambda_0)$ on $\pi^{-1}(\lambda_0)$. Let $t$ be a
transcendental and $K=\bar{\QQ}_p(t)$. Choose a $K$-model
$Z\to\PP^1_K$ of the $m$-cyclic cover of type ${(\bf x;a)}$, where
$x_4=t$. Define $W:=J_Z[p]_\chi(\bar{K})\simeq (\FF_p)^2$. This
defines the {\sl monodromy representation}
\begin{equation}\label{monorep}\rho:\Gal(\bar{K}, K)\to \GL(W)\simeq\GL_2(p).\end{equation}
Let $\Gamma$ be the image of the monodromy representation.  The first
goal of this section is to describe $\Gamma$. These results are due to
V\"olklein \cite{Voelklein93}, \cite{Voelklein95}. We repeat the
proofs as our assumptions are somewhat different.

The action of $\Pi$ on $\pi({\bf a})^{-1}(\lambda_0)$ can be described
via the action of the {\sl Artin braid group} $\B_4$. Recall that
$\B_4$ is generated by the standard generators $Q_1, Q_2, Q_3$ which
act on Nielsen tuples via:
\[ {\bf g}Q_i=(g_1, \ldots, g_{i-1}, g_{i+1},  g_{i+1}^{-1}g_{i} g_{i+1}, g_{i+2}, \ldots, g_r).\]
We follow here the convention of \cite{Voelklein95}. In some other
texts the inverses of the $Q_i$ are used; this does not make much
difference, \cite[p.\ 167]{Voelklein}.  

The {\sl pure Artin braid group} $\B^{(4)}$ is the kernel of the map
$\B_4\to S_4$, sending $Q_i$ to the transposition $(i, i+1)$.  The
pure Artin braid group $\B^{(4)}$ acts naturally on the set
$\E^\ast({\bf a})$. The fundamental group $\Pi$ can be embedded into
the Hurwitz braid group which is a quotient of the Artin braid group.
The action of $\Pi$ on the  Nielsen class $\Niin_4({\bf a})$
is induced from the action of $\B^{(4)}$ on $\E^\ast({\bf a})$. We
refer to \cite{Voelklein95} for details on braid groups.

The pure Artin braid group is generated by the braids
\begin{equation}\label{defbi} b_1=Q_3Q_2Q_1^2Q_2^{-1}Q_3^{-1}, \quad b_2=Q_3Q_2^2Q_3^{-1}, \quad
  b_{3}=Q_3^2.\end{equation} Write $g^h=h^{-1} g h$ and
$[g,h]=g^{-1}h^{-1}gh$.  One checks that
\begin{equation}\label{actionbi}(g_1, g_2, g_3, g_4)b_i=\left\{\begin{array}{ll}
      (g_1^{g_4}, g_2^{[g_1, g_4]}, g_3^{[g_1, g_4]}, g_4^{g_1g_4}) &\mbox{ if }i=1,\\
      (g_1, g_2^{g_4}, g_3^{[g_2, g_4]}, g_4^{g_2g_4}) &\mbox{ if }i=2,\\
      (g_1, g_2, g_3^{g_4}, g_4^{g_3g_4}) &\mbox{ if }i=3.
\end{array}\right.\end{equation}

We start by studying the action of $\B^{(4)}$ on the set
${\E}^\ast({\bf a})$.  The elements of $\E^\ast({\bf a})$
may be represented as
\begin{equation}\label{tuple}([0, \zeta_1], [v_1, \zeta_2], [v_2, \zeta_3], [v_3, \zeta_4]), \quad v_i\in \FF_p,\end{equation}
where $v_3=-(\zeta_2\zeta_3)^{-1}v_1-\zeta_3^{-1}v_2$. Recall that
$\zeta_i=\chi(1)^{a_i}$ and therefore
$\zeta_1\zeta_2\zeta_3\zeta_4=1$. We may identify a tuple
(\ref{tuple}) with a vector ${\bf v}=(v_1, v_2)\in
W:=J(Z)[p]_\chi\simeq (\FF_p)^2$. Here $Z=Z_{\lambda}$ is the
$m$-cyclic cover of type $({\bf x;a})$.  Note that ${\bf v}\in W$
corresponds to a tuple in $\E^\ast({\bf a})$ if and only if ${\bf
  v}\neq 0$.  The action of $\B^{(4)}$ on $\E^\ast({\bf a})$ induces
an action on $W$.

\begin{prop}\label{imrho}
The image in $\GL(W)$ of the action of $\B^{(4)}$ on $W$ is equal to $\Gamma$.
\end{prop}

\proof
\cite[Theorem A]{Voelklein95}
\Endproof

This proposition gives us a concrete way to compute the image of the
monodromy action.

For $i=1,2,3$, define $B_i\in \GL_2(p)$ as the matrix corresponding to the braid $b_i$. As in \cite[Lemma 3]{Voelklein93}, one computes that
\begin{eqnarray}
B_1&=&\left(\begin{array}{cc}\label{B1}
\zeta_2&-1+\zeta_3\\
\zeta_2(-1+\zeta_2)&1-\zeta_2+\zeta_2\zeta_3
\end{array}\right),
\\
B_2&=&\left(\begin{array}{cc}\label{B2}
\zeta_1(1-\zeta_2+\zeta_2\zeta_3)&\zeta_2^{-1}(-1+\zeta_3)(-1+\zeta_1-\zeta_1\zeta_2+\zeta_1\zeta_2\zeta_3) \\
\zeta_1\zeta_2(1-\zeta_2)& 1-\zeta_1+\zeta_1\zeta_2+\zeta_1\zeta_3-\zeta_1\zeta_2\zeta_3
\end{array}\right),\\
B_3&=&\left(\begin{array}{cc}\label{B3}
1&\zeta_1(1-\zeta_3)\\
0&\zeta_1\zeta_2
\end{array}\right).\end{eqnarray}

Let $\bar{\Gamma}$ be the image of $\Gamma$ under the projection
$\GL_2(p)\to \PGL_2(p)$. We show that $\bar{\Gamma}$ is either a group
of order prime-to-$p$ or $ \PSL_2(p)$ or $\PGL_2(p)$. 

\begin{defn}\label{exceptionaldef}
  If $\bar{\Gamma}\neq \PSL_2(p), \PGL_2(p)$, we say that ${\bf a}$ is
  {\sl exceptional}.
\end{defn}

\begin{prop}\label{Hconn}  Suppose that $p>5$. 
\begin{itemize}
\item[(a)] Suppose that $m$ is even and, after permuting the branch
  points, we have ${\bf a}=(a, a, -a+m/2, -a+m/2)$. Then $\Gamma$ is a
  dihedral group of order prime-to-$p$.
\item[(b)] Suppose that ${\bf a}$ is exceptional, but not as in (a).
  The order of $\zeta_i\zeta_j$ in $\FF_p^\times$ is either 2,3,4,5,
  for all $i\neq j$.  If ${\bf a}$ is non-exceptional, ${\Gamma}$
  contains $ \SL_2(p)$.
\item[(c)] Suppose ${\bf a}$ is non-exceptional. The image of $\det:
  \Gamma\to\FF_p^\times$ is generated by $\zeta_1\zeta_2,
  \zeta_1\zeta_3, \zeta_2\zeta_3$.
\item[(d)] Suppose ${\bf a}$ is non-exceptional. The Hurwitz space $H$
  is connected.
\end{itemize}
\end{prop}

\proof 
Part (a) follows easily from the expression for the matrices $B_i$,
cf.\ \cite[Section 1.7]{Voelklein93}.

Suppose that no permutation of ${\bf a}$ is equal to $ (a, a, -a+m/2,
-a+m/2)$.  Recall that a subgroup of $\GL_2(W)$ is called {\sl
  primitive} if it is irreducible and does not permute the summands in
any non-trivial direct sum decomposition of $W$.  It is proved in
\cite[Section 1.7]{Voelklein93} that $\Gamma$ acts primitively on $W$.
We may regard $\Gamma$ as a subgroup of $\PSL_2(p^2)$.  Dickson's
classification of the subgroups of $\PSL_2(p^2)$, \cite[Theorem
8.27]{Huppert}, implies that $\bar{\Gamma}$ is isomorphic to either
$A_4, S_4, A_5, \PSL_2(p)$ or $\PGL_2(p)$.  Therefore, if ${\bf a}$ is
non-exceptional, $\Gamma$ contains $\SL_2(p)$.

Suppose that ${\bf a}$ is exceptional. Then the matrices $B_1, B_2,
B_3$ have order less than or equal to 5. One checks that there exists
vectors $w_i, w_i'\in W$ which are eigenvectors of $B_i$ with
eigenvalue 1 and $(\zeta_i\zeta_4)^{-1}$. If $\zeta_i\zeta_j=1$ for
some $i\neq j$, then $\bar{\Gamma}$ contains an element of order $p>5$
which is impossible. This proves (b).

Suppose that ${\bf a}$ is non-exceptional.  To prove (c), note that the
index of $\SL_2(p)$ in $\Gamma$ is equal to the order of the image of
$\det:\Gamma\to \FF_p^\times$. The subgroup $\det(\Gamma)$ of
$\FF_p^\times$ is generated by the determinant of the matrices $B_i$
for $i=1,2,3$. We have
\[\det(B_i)=\left\{\begin{array}{ll}
\zeta_2\zeta_3&\mbox{ if }i=1,\\
\zeta_1\zeta_3&\mbox{ if }i=2,\\
\zeta_1\zeta_2&\mbox{ if }i=3.
\end{array}\right. \]

Recall that we have identified $\E^*({\bf a})$ with the set of ${\bf
  v}\in W-\{0\}$. Since $\Gamma$ contains $\SL_2(p)$, it follows that
$\Gamma$ acts transitively on $\E^*({\bf a})$. But this implies that
$G$ acts transitively on $\Niin_4({\bf C})=\E^*({\bf a})/(\ZZ/m)$ and
therefore that the Hurwitz space is connected.  \Endproof

It is not clear for which types the group $\bar{\Gamma}$ is either
$A_4,\, S_4$ or $A_5$. The necessary condition of Proposition
\ref{Hconn} is certainly not sufficient. However, for us this is not
important since in the exceptional case $p$ does not divide the order
of the Galois group $G$, since $p>5$.  Therefore, the Galois closure
of $H({\bf a})\to\PP_\lambda^1$ has good reduction to characteristic
$p$.  If $\sum a_i=2m$ we show in Section \ref{grsec} that the Galois
closure of $H({\bf a})\to\PP^1_\lambda$ has bad reduction and we are
not in the exceptional case.

Write $M:=\langle \chi(1)I\rangle$ for the $m$-cyclic subgroup of the
center of $\GL_2(p)$. If ${\bf a}$ is non-exceptional, then
$\SL_2(p)\subset \Gamma\subset\GL_2(p)$. Therefore the index of
$\SL_2(p)$ in $\Gamma$ is equal to the order of $\det(\Gamma)$ . It
follows from Proposition \ref{Hconn}.(c) that $\det(\Gamma)$ is a
subgroup of $\Im(\chi)\simeq \ZZ/m\subset \FF_p^\times$.

\begin{cor}\label{Galpi}
  Suppose that ${\bf a}$ is non-exceptional.  Let $G$ be the Galois
  group of the Galois closure of ${\pi}({\bf a})$. Then
\[G\simeq \left\{ \begin{array}{ll}
    \SL_2(p)& \mbox{ if }\, m\mbox{ is odd},\\
    \PSL_2(p)& \mbox{ if }\, m\mbox{ if even and }\det(\Gamma)\subset\det(M) ,\\
    \PGL_2(p)&\mbox{ if }\, \det(\Gamma)\not\subset\det(M).
\end{array}\right.\] 
\end{cor}

\proof 
Note that the action of $\psi^i\in\ZZ/m$ on $\E^\ast({\bf a})$ corresponds to multiplication by the matrix
\[\left(\begin{array}{cc}
\chi(1)^{i}&0\\
0&\chi(1)^{i}
\end{array}\right)\]
on $W$.  Therefore $G$ is the image of $\Gamma$ under the projection
$\Gamma M\to \Gamma M/M$.

It follows from Proposition \ref{Hconn}(b) that
$\det(\Gamma)\subset\Im(\chi)\simeq \ZZ/m\subset \FF_p^\times$. We
first suppose that $m$ is odd. Then $\ZZ/m$ does not contain an
element of order $2$ and $-I\not\in M$. This implies that
$\det:M\to\FF_p^\times$ is injective and $\det(M)\simeq\ZZ/m$.
Therefore $\det(\Gamma)\subset \det(M)$.  We conclude that the image
$G$ of $\Gamma$ in $\Gamma M\to \Gamma M/M$ is equal to the image of
$\SL_2(p)$ in $\Gamma M\to \Gamma M/M$. Since $-I\not\in \Gamma$, we
have $G\simeq \SL_2(p)$.

Suppose that $m$ is even. Then $-I\in M$. If  $\det(\Gamma)\subset
\det(M)$, then $G$ is isomorphic to
the image of $\SL_2(p)$ in $\Gamma M/M$. Since $-I\in M$, we
conclude that $G\simeq \PSL_2(p)$.

Suppose that $\det(\Gamma)$ is not contained in $\det(M)$. Then $M$
contains $-I$.  The image of $\Gamma$ in $\Gamma M/M$ is strictly
larger than the image of $\SL_2(p)$ in $\Gamma M/M$. We conclude that
$G\simeq \PGL_2(p)$.  \Endproof

\begin{rem} We denote the Galois closure of $\pi({\bf a})$ by
  $\varpi({\bf a}):\HH({\bf a})\to\PP^1$. Following
  \cite{Voelklein95}, we can give a modular interpretation for
  $\HH=\HH({\bf a})$.  We explain the definition and refer to
  \cite{Voelklein95} for details.
  
  Let $\N:=W\rtimes \ZZ/m$, where $\ZZ/m$ acts on $W$ via $\chi$.  Let
  $\H$ be the Hurwitz space over $\bar{\QQ}_p$ parameterizing tuples
  $(f,h,\eta)$, up to isomorphism. Here $f:Y\to \PP^1$ is an
  $\N$-Galois cover of type $({\bf x;a})$ with $x_1=0,\, x_2=1,\,
  x_3=\infty,\, x_4=\lambda$ and $\lambda$ different from the base
  point $\lambda_0$ and $h:\Aut(Y, \PP^1)\to\N$ is an isomorphism.
  Moreover, $\eta$ is an point of $Y$ above $\lambda_0$. We have a
  natural map $\H\to\PP^1_\lambda$, which sends $(f,h,\eta)$ to $\lambda$.

  There is a natural action of $W$ on $\H$. Therefore we can define
  $\H^{(W)}$ as the quotient of $\H$ by $W$. Let $\HH^{(W)}$ be a
  connected component of $\H^{(W)}$. It follows from the results of
  \cite{Voelklein95} that $\varpi^{(W)}:\HH^{(W)}\to\PP^1_\lambda$ is a
  $\Gamma$-Galois cover. The Galois closure $\varpi$ of $\pi$ is the
  quotient of $\varpi^{(W)}$ by $M\cap \Gamma$. We do not need this in
  the rest of the paper.
\end{rem}
  
Our next goal is to describe the ramification of $\pi({\bf a}):H({\bf
  a})\to\PP^1_{\lambda}$. Write $x_1=0, x_2=1, x_3=\infty\in
\PP^1_\lambda$. Let $i\in\{1,2,3\}$.  The points $\pi({\bf
  a})^{-1}(\{x_i\})$ are called the {\sl cusps} of $\pi({\bf a})$. By
\cite[Section 4.2.2]{modular}, there is a 1-1 correspondence between
cusps above $x_i$ for $i\in \{1,2,3\}$ and elements of
\[\Niin_4({\bf a})/\langle b_i\rangle,\]
where the braids $b_i$ defined in (\ref{defbi}) act as in
(\ref{actionbi}).  Define $d_i=m/\gcd(a_i+a_4, m)$ if
$a_i+a_4\not\equiv 0\bmod{m}$ and $d_i=p$ otherwise.

\begin{prop}\label{rampi}
The signature of the Galois closure $\varpi({\bf a}):\HH({\bf a})\to\PP^1$ of $\pi({\bf a})$ is
$(d_1, d_2, d_3)$.
\end{prop}

\proof Choose a cusp $c$ of $H({\bf a})$ above $x_i$ corresponding to
${\bf g}=(g_1, g_2, g_3, g_4)$. For simplicity we suppose that $i=3$,
i.e.\ $x_i=\infty$. The ramification index of the cusp $c$ in
$\pi({\bf a}):H({\bf a})\to\PP^1$ is equal to the length of the orbit
of ${\bf g}$ under $b_3$. It follows from the definition of the action
of $b_i$ that the ramification index of $c$ divides the order of
$g_3g_4$. Let ${\bf g}'=(\psi^{a_1}, \psi^{a_2}, \phi\psi^{a_3},
g_4')$, where $g_4'$ is chosen such that the product of the $g_i'$ is
one.

We distinguish two cases. First suppose that $a_3+a_4=m$. Then for
every Nielsen tuple ${\bf g}$, the order of $g_3g_4$ divides $p$. One
easily checks that the orbit of ${\bf g}'$ has length $p$. This
implies that the ramification index of $\infty$ in $\varpi({\bf a})$ is $p$.

Now suppose that $a_3+a_4<m$. Then the order of $g_3g_4$ divides
$d_3m/\gcd(a_3+a_4, m)$. Since the orbit of ${\bf g}'$ has length
$d_3$, we conclude that $\infty$ is ramified in $\varpi({\bf a})$ of
order $d_3$.  \Endproof

\section{A good reduction theorem}
\label{grsec}
Let ${\bf a}=(a_1, a_2, a_3, a_4)$ be a type. To  ${\bf a}$ we
associated in Section \ref{Hurwitz0} a Galois cover $\varpi({\bf
  a}):\HH({\bf a})\to\PP^1$. In this section we determine when
$\varpi({\bf a})$ has good reduction.

Let $K=\bar{\QQ}_p(t)$, where $t$ is a transcendental element. Denote
by $G_K$ the absolute Galois group of $K$. Recall from Proposition
\ref{imrho} that the monodromy representation $\rho:G_K\to\GL_2(p)$
has image $\Gamma$. Let $L'/K$ be the corresponding $\Gamma$-Galois
extension. Recall that the Galois group $G$ of $\varpi({\bf a})$ is
the quotient of $\Gamma$ by the (central) subgroup $\Gamma\cap M$ of
$\Gamma$. Let $L$ be the sub field of $L$ of invariants under
$\Gamma\cap M$.  Denote the Gau\ss-valuation of K by $\nu_0$; this
valuation corresponds to the original component of the stable model of
$\varpi({\bf a})$.

\begin{lem}\label{grcrit}
The cover $\varpi({\bf a})$ has good reduction if and only if $\nu_0$ is
unramified in $L/K$.
\end{lem}

\proof The forward implication is obvious. Suppose that $\nu_0$ is
unramified in $L/K$. Then the restriction of the stable model of
$\varpi({\bf a})$ to the original component is separable. It is well
known that this implies that $\varpi({\bf a})$ has good reduction,
\cite[Proposition 1.1.4]{modular}.  \Endproof

\begin{thm} \label{grthm}
The cover $\varpi({\bf a})$ has good reduction if and only if
$a_1+a_2+a_3+a_4\neq 2m$.
\end{thm}

\proof Let $R$ be the ring of integers of the completion of $K$ at
$\nu_0$. Let $k$ be the residue field of $R$. Denote by $Z_{0,R}$ the
$m$-cyclic cover of $\PP^1_R$ of type ${\bf a}$ branched at
$0,\,1,\,\infty,\,t$. Let $\Zb_0$ be the special fiber of $Z_{0,R}$.
There is an exact sequence of group schemes
\begin{equation}\label{splitseq}1\to J_{\Zb_0}[p]_\chi^\loc\to J_{\Zb_0}[p]_\chi\to J_{\Zb_0}[p]_\chi^\et\to
  1.\end{equation}

Suppose that $\sum_i a_i=m$, i.e.\ we are in the multiplicative case
of Definition \ref{casedef}. It follows from Proposition \ref{mindim}
that $ J_{Z_{0,R}}[p]_\chi$ is \'etale. This implies that the action of
$\Gamma$ on $ J_{Z_{0,R}}[p]_\chi$ is unramified, since $R$ is a
Henselian local ring. Lemma \ref{grcrit} implies therefore that
$\varpi({\bf a})$ has good reduction.

Suppose that $\sum_i a_i=3m$, i.e.\ we are in the \'etale case. The
group scheme $ J_{Z_{0,R}}[p]_\chi$ is dual to
$ J_{Z_{0,R}}[p]_{\chi^{-1}}$. But $ J_{Z_{0,R}}[p]_{\chi^{-1}}\simeq
 J_{Z^\ast_{0,R}}[p]_\chi$, where $Z^\ast_{0,R}$ is the $m$-cyclic
cover of type ${\bf a}^\ast=(m-a_1, m-a_2, m-a_3, m-a_4)$ branched at
$0,\,1,\,\infty,\, t$. Since $\sum_i(m-a_i)=m$, we may apply the
theorem in the multiplicative case to ${\bf a}^\ast$. The theorem in
the \'etale case follows by duality.

Suppose that $\sum_ia_i=2m$, i.e.\ we are in the mixed case. Since the
fourth branch point $x_4=t$ is generic, it follows from Proposition
\ref{genord} that both $ J_{\Zb_0}[p]_\chi^\loc$ and
$ J_{\Zb_0}[p]_\chi^\et$ have rank $p$. Choose a valuation $\nu'$ of
$L'$ above $\nu_0$ and write $\Gamma_{\nu'}$ for its decomposition
group. Suppose that $\varpi({\bf a})$ has good reduction. The
$\Gamma$-Galois cover $\varpi'({\bf a})$ corresponding to $L'/K$ has
good reduction also. (Here we use that the index of $G$ in $\Gamma$ is
prime-to-$p$.) We conclude that $\Gamma$ is equal to $\Gamma_{\nu'}$
and that $\Gamma$ acts faithfully on the stable reduction of
$\varpi'({\bf a})$. In particular, $\Gamma$ acts transitively on the
subgroups of $ J_{\Zb_0}[p]_\chi$. But the action of $\Gamma$ on
$ J_{\Zb_0}[p]_\chi$ respects the splitting (\ref{splitseq}) of
$ J_{\Zb_0}[p]_\chi$ into its local and \'etale parts. This gives a
contradiction.
\Endproof

\section{The ${\bf a}$-Hasse invariant}\label{Hassesec}
In Sections \ref{Hassesec} and \ref{defodatumsec}, we suppose that
$\sum_i a_i=2m$, i.e.\ $\varpi({\bf a})$ has bad reduction. Our goal
is to compute the stable reduction of $\varpi({\bf a})$. In this
section we define the ${\bf a}$-Hasse invariant $\Phi_{\bf a}$. This
is a generalization if the classical Hasse invariant for elliptic
curves, cf.\ \cite[Section 12.4]{KatzMazur}. The ${\bf a}$-Hasse
invariant satisfies a hypergeometric differential equation, just as
the classical Hasse invariant. In Section \ref{defodatumsec} we show
that $\Phi_{\bf a}$ essentially determines the stable reduction.

Let ${\bf a}$ be a type such that $a_1+a_2+a_3+a_4=2m$. Let $\Zb_0$ be the
$m$-cyclic cover of $\PP^1_k$ branched at $0,\,1,\,\infty,\, \lambda$ of
type ${\bf a}$. Since we are in the mixed case, Lemma \ref{dimLi}
  implies that the $\chi$-isotypical part $H^0(\Zb_0, \Omega)_\chi$
  has $k$-dimension one. Recall that $\C:H^0(\Zb_0,\Omega)_\chi\to
  H^0(\Zb_0,\Omega)_\chi$ is a bijection for all but finitely many
  $\lambda$, Proposition \ref{genord}.

\begin{defn}\label{ssdef}
  The $\bar{\lambda}$ for which $\C:H^0(\Zb_0,\Omega)_\chi\to
  H^0(\Zb_0,\Omega)_\chi$ is a bijection are called ${\bf a}$-{\sl
    ordinary}. If $\bar{\lambda}$ is not ${\bf a}$-ordinary it is
  called ${\bf a}$-{\sl supersingular}.  Write $\Lambda({\bf a})$ for
  the set of ${\bf a}$-supersingular $\lambda$-values.
\end{defn}

For $m=2$, the only possible type is $(1,1,1,1)$. In this case $\Zb_0$
is an elliptic curve and $\lambda$ is ${\bf a}$-supersingular if and
only if $\Zb_0$ is supersingular.

Now suppose that $\lambda=t$ corresponds to the generic point of
$\PP^1_k$. Then $\Zb_0$ is given by an equation
$z^m=x^{a_1}(x-1)^{a_2}(x-\lambda)^{a_4}$. Define
\[\omega=\frac{z\, {\rm d}x}{x(x-1)(x-\lambda)}\in H^0(\Zb_0, \Omega^1)_\chi.\]
Lemma \ref{dimLi} implies that $\dim_k H^0(\Zb_0,
\Omega^1)_\chi=\dim_k H^1(\Zb_0, {\cal O}_{\Zb_0})_{\chi^{-1}}=1$.
Therefore $\omega$ is a basis of $H^0(\Zb_0, \Omega^1)_\chi$.

\begin{defn}\label{Phidef}
Suppose that $a_1+a_2+a_3+a_4=2m$. Define the ${\bf a}$-Hasse invariant as
\[\C\omega=\Phi_{\bf a}(\lambda)^{(1/p)}\omega.\]
\end{defn}

\begin{lem}\label{phi}
  Write $a^\ast_i=m-a_i$ and $\alpha=(p-1)/m$. Then
\begin{equation}\label{formulaphi}\Phi_{\bf a}(\lambda)=(-1)^{\alpha a_3}\sum_{i+j=\alpha a_3} \left(\begin{array}{c}\alpha a^\ast_2 \\ i\end{array}
\right)\left(\begin{array}{c}\alpha a^\ast_4 \\ j\end{array}\right) \lambda^j.\end{equation}
\end{lem}

\proof
Write
\[\omega=\frac{z^p}{x^{p}(x-1)^p(x-\lambda)^p}x^{p-\alpha a_1}(x-1)^{p-1-\alpha a_2}(x-\lambda)^{p-1-\alpha a_4} \, \frac{{\rm d} x}{x}\]
and 
\[e:=x^{p-\alpha a_1}(x-1)^{p-1-\alpha a_2}(x-\lambda)^{p-1-\alpha a_4}=x^{1+\alpha a^\ast_1}(x-1)^{\alpha a^\ast_2}(x-\lambda)^{\alpha a^\ast_4}=\sum_i e_i x^i.\]
It follows from standard properties of the Cartier operator that
\[\C\omega=\frac{z}{x(x-1)(x-\lambda)}\left( \sum_i e_{pi}^{1/p}
  x^i\right)\frac{{\rm d} x}{x},\] \cite{SerreMex}.  Note that
$\deg{e}=1+\alpha(a^\ast_1+a^\ast_2+a^\ast_4)<1+\alpha(a^\ast_1+a^\ast_2+a^\ast_3+a^\ast_4)=1+2\alpha
m=2p-1$ and $\deg(e)\geq
1+\alpha(a^\ast_1+a^\ast_2+a^\ast_3+a^\ast_4)-\alpha m=1+\alpha m=p$.
Therefore
\[\C\omega=\frac{z(e_p^{1/p}x)}{x(x-1)(x-\lambda)}\frac{{\rm d} x}{x}=e_p^{1/p}\omega.\]
One easily checks that $e_p$ is equal to the right hand side of
(\ref{formulaphi}).  \Endproof

We find back the classical Hasse invariant $\Phi$ for $m=2$ and ${\bf
  a}=(1,1,1,1)$. It is well known that $\Phi$ is the solution to a
hypergeometric differential equation. The same holds for the ${\bf
  a}$-Hasse invariant.

\begin{prop}\label{hgde}
  The polynomial $\Phi_{\bf a}(\lambda)$ is a solution to the
  hypergeometric differential equation
\begin{equation}\label{hgdeeq}\lambda (1-\lambda)u''+[C-(A+B+1)\lambda]u'-ABu=0,\end{equation}
where $A\equiv -\alpha a_3 \bmod{p}$ and  $B\equiv \alpha (a_4-m)
\bmod{p}$ and $C\equiv -\alpha(a_2+a_3) \bmod{p}.$
\end{prop}

\proof Write $\Phi_{\bf a}(\lambda)=\sum_{j} \varphi_j \lambda^j.$ To
prove the proposition, we have to show that the coefficients
$\varphi_j$ satisfy the recursion relation
$\varphi_{j+1}/\varphi_j=(A+j)(B+j)/(C+j)(1+j)$, \cite[Section
III.2]{Yoshida}. Here $j$ and $j+1$ are in the range where the
coefficients are not congruent to zero.  Obviously,
\[\frac{\varphi_{j+1}}{\varphi_j}=\frac{(-\alpha a_3+j)(-\alpha a^\ast_4+j)}{(-\alpha (a_2+a_3)+j)(j+1)}\equiv\frac{(A+j)(B+j)}{(C+j)(1+j)} \pmod{p},\]
with $A,B,C$ as in the statement of the proposition.  \Endproof

\begin{cor}\label{zerosphi}
\begin{itemize}
\item[(a)] The zeros of $\Phi_{\bf a}(\lambda)$
  different from $0$ and 1 are simple. 
\item[(b)] At $\lambda=0$, the polynomial $\Phi_{\bf a}(\lambda)$ has
  a zero of order $\max(\alpha(a_2+a_3-m), 0)$. 
\item[(c)] At $\lambda=1$, the polynomial $\Phi_{\bf a}(\lambda)$ has
  a zero of order $\max(\alpha(a_1+a_3-m), 0)$. 
\end{itemize}
\end{cor}

\proof Part (a) follows immediately from Proposition \ref{hgde} as the
hypergeometric differential equation (\ref{hgdeeq}) has essential
singularities only at $0,1,\infty$. Part (b) is obvious. Part (c)
follows from Part(b) by interchanging the branch points $x_1$ and $x_2$.
\Endproof

\begin{cor}\label{numss}
 The number of ${\bf a}$-supersingular $\lambda$-values is
\[\min(\alpha a^\ast_4, \alpha a_3)-\max(\alpha(a_2+a_3-m), 0)-\max(\alpha (a_1+a_3-m), 0).\] 
\end{cor}

\proof The degree of $\Phi_{\bf a}$ is equal to $\min(\alpha a^\ast_4,
\alpha a_3)$. Therefore the corollary follows from Corollary
\ref{zerosphi}. One easily checks from the formula that the number of
${\bf a}$-supersingular $\lambda$-values is non-zero. \Endproof

\section{Computation of the deformation datum}\label{defodatumsec}
The goal of this section is to describe the stable reduction of
$\varpi({\bf a}):\HH({\bf a})\to\PP^1$ By Theorem \ref{specialthm},
this amounts to computing the deformation datum associated to
$\varpi({\bf a})$. 

The notation is as in Section \ref{grsec}. In particular, we suppose
that $\sum_i a_i=2m$.  Let $k$ be an algebraically closed field of
characteristic $p$. We denote by $P$ the Sylow $p$-subgroup of
$\Gamma$ consisting of upper triangular matrices with ones on the
diagonal.  Choose a valuation $\nu$ of $L'$ above $\nu_0$ such that
the inertia group $I_\nu$ of $\nu$ is $P$. This is possible as
$\varpi({\bf a})$ has bad reduction by Theorem \ref{grthm}. Let
$\Gamma_\nu$ be the decomposition group of $\nu$. By our assumption on
$\nu$, the group $\Gamma_\nu$ is a subgroup of the Borel group $\B$ of
$\GL_2(p)$ consisting of upper triangular matrices. Define
\[\Gamma_\nu^\res=\Gamma_\nu/I_\nu.\]
The monodromy representation (\ref{monorep}) induces a representation
\[\rho^\res:\Gamma_\nu^\res\to\B/P=\FF_p^\times\times\FF_p^\times.\]
We may identify $\rho^\res$ with two characters
$\xi_i:\Gamma_\nu^\res\to\FF_p^\times$, where $\xi_i$ is the
composition of $\rho^\res$ with the $i$th projection.

\begin{prop}\label{monoprop}
  Denote by $G_{k(t)}$ the absolute Galois group of $k(t)$.  Under the
  natural identification
\begin{eqnarray*}
  k(t)^\times/(k(t)^\times)^{p-1}&\simeq& \Hom(G_{k(t)}, {\boldsymbol \mu}_{p-1})\\
  u&\mapsto& \left[\sigma\mapsto \ 
    \frac{^\sigma\!(\sqrt[p-1]{u})}{\sqrt[p-1]{u}}\right]
\end{eqnarray*}
the character $\xi_2$ (resp.\ $\xi_1$) corresponds to $\Phi_{{\bf
a}}(t)$ (resp.\ $\Phi_{\bf a^\ast}(t)^{-1}$).  Here $\Phi_{\bf a}$ is
the ${\bf a}$-Hasse invariant and ${\bf a}^\ast:=(m-a_1, m-a_2, m-a_3,
m-a_4)$ is the dual type.
\end{prop}

\proof Let ${\bf x}=(0,1,\infty, t)$ and write $\Zb_0$ for the
$m$-cyclic cover of type ${\bf (x;a)}$ over $k(t)$.  It follows from
(\ref{splitseq}) that the character $\xi_2$ (resp.\ $\xi_1$)
corresponds to the action of $G_{k(t)}$ on $J_{\Zb_0}[p]_\chi^{\et}$
(resp.\ $J_{\Zb_0}[p]_\chi^{\loc}$).

We start by computing $\xi_2$. For this we use the canonical
identification
\[J_{\Zb_0}[p]_\chi^{\et}\simeq H^0(\Zb_0, \Omega)_\chi^{\C},\]
as in Section \ref{duality}.  The
differential
\[\omega:=\frac{z\, {\rm d} x}{x(x-1)(x-t)}\]
is a basis of $H^0(\Zb_0, \Omega)_\chi$, cf.\ Section \ref{Hassesec}.  It
follows from Definition \ref{Phidef} that $\C\omega=\Phi_{\bf
  a}(t)^{1/p}\omega.$ Since $\C c\omega=c^{1/p}\C\omega$, we see that
\[\C \Phi_{\bf a}(t)^{\frac{1}{p-1}} \omega=\Phi_{\bf a}(t)^{\frac{1}{p-1}} \omega.\]
We conclude that $\xi_2$ corresponds to the $(p-1)$-cyclic extension
of $k(t)$ obtained by adjoining a $(p-1)$th root of the ${\bf
a}$-Hasse invariant $\Phi_{\bf a}(t)$.

To prove the statement of the proposition for $\xi_1$, we remark that
$J_Z[p]_\chi^{\loc}$ is {dual} to $ J_Z[p]_{\chi^{-1}}^{\et}$.
Therefore the above argument applied to the character $\chi^{-1}$
shows that $\xi_1$ corresponds to the inverse of the Hasse invariant
$\Phi_{{\bf a}^\ast}(t)$ corresponding to the dual type ${\bf
  a}^\ast=(m-a_1, m-a_2, m-a_3, m-a_4)$.  \Endproof

Let $\bar{W}_0$ be the original component of the stable reduction of
the $G$-Galois cover $\varpi({\bf a}):\HH({\bf a})\to\PP^1_\lambda$.
We choose a component $\bar{U}_0$ of $\bar{\HH}({\bf a})$ above
$\bar{W}_0$.  The map $\bar{U}_0\to\bar{W}_0$ factors through a curve
$\bar{V}_0$ such that $\bar{U}_0\to\bar{V}_0$ is a $\mup$-torsor and
$\bar{V}_0\to\bar{W}_0$ is a tame cover, cf.\ Section \ref{3ptsec}.
The following theorem describes the cover $\bar{V}_0\to\bar{W}_0$.

\begin{thm}\label{defodatumA}
  Suppose that  ${\bf a}$ is non-exceptional, Definition
  \ref{exceptionaldef}. Write $d:=\gcd(m,
  a_1+a_3, a_2+a_3)$. If $m\neq 2$, we  assume that $d\neq m$. This is
  no restriction.
\begin{itemize}
\item[(a)] Suppose that $m=2$. Then $G=\PSL_2(p)$ and the cover
  $\bar{V}_0\to\bar{W}_0$ is given by
\[\theta^{(p-1)/2}=\Phi_{\bf a}.\] 
\item[(b)] Suppose that $G=\SL_2(p)$. Write $m/d=1+2j$.  The cover
  $\bar{V}_0\to\bar{W}_0$ is given by
\[\theta^{p-1}=\Phi_{\bf a}^{1+j}\Phi_{\bf a^\ast}^{-j}.\]
\item[(c)] Suppose $G=\PSL_2(p)$ or $G=\PGL_2(p)$ and $d$ is even.
  Then there exists a function $g\in k(t)$ such that
\[\Phi_{\bf a}^{1+m/d}(t)\Phi_{\bf
  a^\ast}^{1-m/d}(t)=g^2.
\]
The cover $\bar{V}_0\to\bar{W}_0$ is given by
\[\theta^{(p-1)/2}=g.\]
\item[(d)] Suppose that $G=\PGL_2(p)$ and $d$ is odd.  The cover
  $\bar{V}_0\to\bar{W}_0$ is given by
\[\theta^{p-1}=\Phi_{\bf a}^{1+m/(2d)}\Phi_{\bf a^\ast}^{1-m/(2d)}.\]
\end{itemize}
\end {thm}

\proof
We first note  that
\[\Phi_{{\bf a}^\ast}(t)=(-1)^{\alpha a_3^\ast}\sum_{i+j=\alpha
  a_3^\ast} \left(\begin{array}{c}{\alpha a_2}\\{i}\end{array}\right)\left(\begin{array}{c}{\alpha a_4}\\{j}\end{array}\right)\lambda^j.\]
The polynomial $\Phi_{{\bf a}^\ast}(t)$ has a zero at $t=0$ of order
$\max(\alpha(m-a_2-a_3), 0)$ and a zero at $t=1$ of order
$\max(\alpha(m-a_1-a_3), 0)$.

Since the group scheme $J_{\Zb_0}[p]_\chi$ is dual to
$J_{\Zb_0}[p]_{\chi^{-1}}$, it is obvious that $J_{\Zb_0}[p]_\chi$ is
local-local if and only if $J_{\Zb_0}[p]_{\chi^{-1}}$ is local-local.
This implies that $t\neq 0,1$ is a zero of $\Phi_{{\bf a}}(t)$ if and
only if it is a zero of $\Phi_{{\bf a}^\ast}(t)$. Recall from
Corollary \ref{zerosphi} that all zeros
 of $\Phi_{\bf a}$ and $\Phi_{\bf a^\ast}$, except 0 and 1, are simple.
Therefore
\[\frac{\Phi_{{\bf a}}(t)}{\Phi_{{\bf a}^\ast}(t)}=t^{\pm \alpha(a_2+a_3-m)}(t-1)^{\pm \alpha (a_1+a_3-m)}.\]

\bigskip\noindent{\bf Case (a)}. Suppose that $m=2$. The only
possibility for the type is ${\bf a}=(1,1,1,1)$. Since $a_i+a_j=m$,
for all $i$ and $j$, we conclude that $\Phi_{\bf a}=\Phi_{\bf
  a}^\ast$. This implies the theorem in this case. This is well
known, see for example \cite{KatzMazur}.

Suppose now that $m\neq 2$. Then the exist $i\neq j$ such that
$a_i+a_j\neq m$. Therefore, after permuting the branch points if
necessary, we may assume that $d\neq m$. Let $L_2$ (resp.\ $L_1$) be
the extension of $k(t)$ corresponding to adjoining a $(p-1)$th root
$\theta_i$ of $\Phi_{{\bf a}}(t)$ (resp.\ $1/\Phi_{\bf a^\ast}(t)$).
Let $\tilde{L}=L_1\otimes_{k(t)} L_2$.  Choose a primitive $(p-1)$th
root of unity $\xi\in k$. Choose generators $\sigma_i$ for the Galois
group of $\tilde{F}$ over $k(t)$ such that
\[\sigma_1 \theta_1=\xi \theta_1,\qquad
\sigma_1\theta_2=\theta_2,\qquad
\sigma_2\theta_1=\theta_1,\qquad\sigma_2\theta_2=\xi \theta_2.\] 
For $\ell|(p-1)$, denote by $L_i^{(\ell)}$ the
(unique) subfield of $L_i$ containing $k(t)$ whose degree over $k(t)$
is $\ell$. Note that $L_1^{(\ell)}\simeq L_2^{(\ell)}$ if and only if
$\Phi_{{\bf a}}\Phi_{{\bf a}^\ast}^i=\delta^\ell$ for some $\delta\in k(t)$ and
some $i$ prime to $(p-1)$.  Let $d$ be as in the statement of the
theorem and $\alpha=(p-1)/m$. Then $[L_1\cap L_2:k(t)]=\alpha d$.
Therefore the algebra $\tilde{L}$ is the product of $\alpha d $ fields
$L(\eta)$ with
\[L(\eta)=k(t)[\theta_1, \theta_2|\, \theta_2^{p-1}=\Phi_{\bf a}, \quad 
\theta_1^{p-1}=1/\Phi_{\bf a^\ast}, \quad
(\theta_1\theta_2)^{m/d}=\eta \delta],\] where $\eta$ is a $\alpha d$th
root of unity and $\delta^{\alpha d}={\Phi_{\bf a}}/{\Phi_{\bf a^\ast}}.$
The group $G(\eta)=\Gal(L(\eta), k(t))$ is generated by
$\sigma_1\sigma_2^{-1}$ and $\sigma_1^{\alpha d}$.

Recall that $M$ is the subgroup of $\GL_2(p)$ consisting of the scalar
matrices $xI$ with $x\in \ZZ/m\subset \FF_p^\times$.  To find the
subfield of $L(\eta)$ which corresponds to $\bar{V}_0$, we have to
take invariants under the restriction of $M\cap \Gamma$ to $L(\eta)$,
cf.\ Section \ref{Hurwitz0}. From now on, we suppose $\eta=1$. This is no
restriction.

Define $J$ as the intersection of $M\cap \Gamma$ with $G(\eta)$.  The
group $M\cap \Gamma$ is generated by $(\sigma_1\sigma_2)^\alpha$.
Therefore $I$ is generated by $(\sigma_1\sigma_2)^{\alpha d/2}$ if $d$
is even and by $(\sigma_1\sigma_2)^{\alpha d}$ if $d$ is odd.

\bigskip\noindent {\bf Case (b)} Suppose that $G=\SL_2(p)$. Corollary
\ref{Galpi} implies that $m$ is odd. Write $m/d=1+2j$.  The degree of
$L(\eta)^J$ over $k(t)$ is $p-1$.  Define
$\theta:=\theta_1^{j}\theta_2^{1+j}$. We claim that $\theta$ generates
$L(\eta)^J$ over $k(t)$. It is clear that $\theta\in L(\eta)^J$. Since
we assumed that $m\neq d$, the integer $j$ does not divide $p-1$.
Therefore $\theta^i\not\in k(t)$, for every $0<i<p-1$ and $\theta$ generates $L(\eta)^J$ over $k(t)$.  We have
\begin{equation}\label{eqa}\theta^{p-1}=\Phi_{\bf a}^{1+j}\Phi_{\bf
  a^\ast}^{-j}=x^{b_1}(x-1)^{b_2}\prod_{\lambda\in\Lambda({\bf a})}(x-\lambda).\end{equation}
The $b_i$ are expressions in terms of the $a_j$ which we leave to the
reader to compute.

\bigskip\noindent {\bf Case (c)}. Suppose that $G=\PSL_2(p)$.
Corollary \ref{Galpi} implies that $m$ is even and that $\det(\Gamma)$
is a subgroup of $\det(M)$. Therefore the order of $\det(\Gamma)$
divides $m/2$. This implies that $\gcd(m, a_1+a_3, a_1+a_3, a_3+a_4)$
is even and  so  $d$ is even. We conclude that $J$ is generated by
$(\sigma_1\sigma_2)^{\alpha d/2}$ and that the degree of $L(\eta)^J$
over $k(t)$ is $(p-1)/2$.

One checks that $L(\eta)^J$ is generated over $k(t)$  by $\theta:=\theta_1^{-1-m/d}\theta_2^{+1-m/d}$ which satisfies
\[\theta^{p-1}=\Phi_{\bf a}^{1+m/d}\Phi_{\bf
  a^\ast}^{1-m/d}=x^{b_1}(x-1)^{b_2}\prod_{\lambda\in\Lambda({\bf
    a})}(x-\lambda)^2.\]
 Then
 \[b_i=\left\{\begin{array}{cc}
\alpha(a_i+a_4-m)(1-\frac{m}{d})& \mbox{ if } a_i+a_4>m,\\
\alpha(m-a_i-a_4)(1+\frac{m}{d})& \mbox{ if } a_{i}+a_4<m.
\end{array}\right.\]
In case $a_{i}+a_4=m$ the branch point $x_i$ is unramified
in $\bar{V}_0\to\bar{W}_0$.

Since $d$ is even, also $b_1$ and $b_2$ are even. We conclude that the
$(p-1)/2$-cyclic cover $\bar{V}_0\to \bar{W}_0$ is given by the
equation
\begin{equation}\label{eqb}(\theta)^{{(p-1)}/{2}}=x^{{b_0}/{2}}(x-1)^{{b_1}/{2}}\prod_\lambda
(x-\lambda).\end{equation}

We leave it for the reader to verify that if $G=\PGL_2(p)$ and $d$ is
even the cover $\bar{V}_0\to \bar{W}_0$ is given by the same equation.

\bigskip\noindent {\bf Case (d)} Suppose that $G=\PGL_2(p)$ and that
$d$ is odd. It follows from the proof of Corollary \ref{Galpi} that in
this case $m$ is even. Define $j=m/(2d)$. Since $d$ is odd, we know
that $J$ is generated by $(\sigma_1\sigma_2)^{\alpha d}$ and that the
degree of $L(\eta)^J$ over $k(t)$ is $p-1$. One checks that
$L(\eta)^J/k(t)$ is generated by
$\theta:=\theta_1^{-1+j}\theta_2^{1+j}$ which satisfies
\begin{equation}\label{eqc}\theta^{p-1}=\Phi_{\bf a}^{1+j}\Phi_{\bf
    a^\ast}^{1-j}=x^{b_1}(x-1)^{b_2}\prod_{\lambda\in\Lambda({\bf
      a})}(x-\lambda)^2.\end{equation} \Endproof

Recall the following from Section \ref{stablesec}. The stable
reduction $\bar{\varpi}({\bf a}):\bar{\HH}({\bf a})\to \bar{\PP}$ of
$\varpi({\bf a})$ has a very simple structure. The curve $\bar{\PP}$
is a comb. It consist of tails $W_b$ for $b\in\BB$ which intersect the
original component in one point which we denote by $\tau_b$. There are
two types of tails: new tails and primitive tails. A tail $W_b$ is
primitive if $W_b$ contains the specialization of one of the branch
points of $\varpi({\bf a})$ whose ramification index is prime-to-$p$.
Therefore $x_i$ for $i\in\{1,2,3\}$ specializes to a primitive tail if
$a_i+a_4\neq m$, Proposition \ref{rampi}. All other tails are called
the new tails.  They intersect the original component $\bar{W}_0$ in an
${\bf a}$-supersingular $\lambda\in\Lambda({\bf a})$.  We may regard
the set $\BB$ indexing the tails as a subset of $\{1,2,3\}\cup
\Lambda({\bf a})$.

 To a tail $W_b$ of $\bar{\PP}$ we associated a ramification
invariant $\sigma_b=h_b/m_b$ which describes the ramification above
the intersection point of $W_b$ with $\bar{W}_0$.  As a consequence of
Theorem \ref{defodatumA}, we can now describe the ramification
invariants of the tails of $\bar{\PP}$.

Recall that we have chosen an irreducible component $\bar{U}_0$ of
$\bar{\HH}({\bf a})$ such that the inertia group of $\bar{U}_0$ is the
Sylow $p$-subgroup $P$ of $G$. The map $\bar{U}_0\to \bar{W}_0$
factors through $\bar{V}_0$ and $\bar{U}_0\to \bar{V}_0$ is a
$\mup$-torsor.  Write $\bar{U}_0'$ for the quotient of $\bar{U}_0$ by
the prime-to-$p$ part of the center of its decomposition group and
write $\bar{V}_0'$ for the quotient of $\bar{U}_0'$ by $\mup$.  Then
$\bar{V}_0'\to\bar{W}_0$ is a cyclic cover, which is branched at $0,
1, \infty, \lambda$ for $\lambda\in\Lambda({\bf a})$. Let $n$ be the
order of this cover.

Remark \ref{typehurw} implies that $0<\sigma_i<1$ for $i\in\{1,2,3\}$
and $1<\sigma_\lambda<2$ for $\lambda\in\Lambda({\bf a})$. Moreover,
the ramification invariants are related to the type of the cover
$g_0':\bar{V}_0'\to\bar{W}_0$, \cite{Stefan01}. The relation is as
follows. Write $\sigma_b=h_b/n_b$, where $n_b$ is the ramification
index of $\tau_b$ in $g_0'$. Let $(\beta_b, \tau_b)$ be the type of
$g_0'$, as defined in Definition \ref{deftype}. Then $h_b n/n_b\equiv
\beta_b\bmod{n}$. The invariants $\beta_b$ are easily computed from
the explicit equation for $g_0'$.  This determines the ramification
invariants.

We will make this explicit in case $\lambda\in\Lambda({\bf a})$.
Suppose first that we are in Case (b) of Theorem \ref{defodatumA},
i.e.\ $G=\SL_2(p)$.  Then $\bar{V}_0\to\bar{W}_0$ has order $p-1$. The
normalizer $N_G(P)$ of the Sylow $p$-subgroup $P$ in $\SL_2(p)$ has
order $p(p-1)$, therefore the decomposition group of $\bar{U}_0$ is
the full normalizer $N_G(P)$. The order of the prime-to-$p$
centralizer of $N_G(P)$ is two. Therefore the degree of
$\bar{V}_0'\to\bar{W}_0$ is $(p-1)/2$. It follows from (\ref{eqa})
that the integer $\beta_\lambda=1$, for every $\lambda\in\Lambda({\bf
  a})$.  The ramification invariant is now determined by the fact that
$(p-1)/2=n_{\lambda}<h_\lambda<2n_\lambda=(p-1)$ and $2h_\lambda\equiv
1\bmod{(p-1)}$. We conclude that $\sigma_\lambda=(p+1)/(p-1)$.  One
checks that the same holds in the other cases.

We note that the fact that all new tails have ramification invariant
$(p+1)/(p-1)$ is a consequence of the fact that the ${\bf a}$-Hasse
invariant is a solution of a hypergeometric differential equation.
The decomposition group of a new tail is either $\PSL_2(p)$ or
$\SL_2(p)$.  One can show the following. Let $f$ be a Galois cover of
$\PP_k^1$ with Galois group either $\PSL_2(p)$ or $\SL_2(p)$.  Suppose
that $f$ is branched only at $\infty$ and that the ramification
invariant $\sigma$ of $f$ satisfies $1<\sigma<2$. Then
$\sigma=(p+1)/(p-1)$.  The proof uses ideas from \cite[Section
3.3]{RRR}.

\begin{cor}\label{defodatumB}
The differential $\omega$ is given by
\begin{equation}\label{omegahurw}\omega=u \frac{\theta \,{\rm d}x}{x(x-1)},\end{equation}
where $\theta$ and $x$ are as in Theorem \ref{defodatumA} and $u\in k^{\times}$.
\end{cor}

\proof 
This is immediate.
 \Endproof
 
 An important property of $\omega$ is that $\C\omega=\omega$. This
 determines the constant $u$ in (\ref{omegahurw}), up to an element of
 $\FF_p^{\times}$. (Actually, $u\in\FF_p^\times$.)  One can also check
 directly that the differential defined by (\ref{omegahurw}) with
 $u\in\FF_p^\times$ is fixed by the Cartier operator, for example
 by using the method of \cite[Section 4.3]{Ihara74}. It turns out that
 $\C\omega=\omega$ is equivalent to the fact that $\Phi_{\bf a}$
 satisfies the hypergeometric differential equation of Proposition
 \ref{hgde}.

\section{Reduction of the Hurwitz spaces}\label{Hurwitzmixed}
In the previous sections, we computed the stable reduction of
$\varpi({\bf a}):\HH({\bf a})\to\PP^1_\lambda$ as Galois cover. We did
not give an interpretation of the stable model $\bar{\HH}({\bf a})$ as
moduli space; it is not clear whether the model $\bar{\HH}({\bf a})$
has a ``reasonable'' interpretation as a moduli space. We merely used
the interpretation of $\HH({\bf a})$ (or rather of its quotient
$H({\bf a})$) as Hurwitz space to determine the stable reduction of
the cover $\varpi({\bf a}):\HH({\bf a})\to\PP^1_\lambda$.

Let $H({\bf a})/K$ be the Hurwitz space parameterizing metacyclic
covers of type ${\bf (x;a)}$, as defined in Section \ref{Hurwitz0}.
Here $K$ is some finite extension of $\QQ_p$. Write $\pi({\bf
  a}):H({\bf a})\to\PP^1_\lambda$ for the projection to the
$\lambda$-line. In \cite{modular}, we find a definition of a
compactification $\bar{H}({\bf a})$ of $H({\bf a})$ over the ring of
integers $R$ of $K$. Roughly speaking, $\bar{H}({\bf a})\otimes_R
\FF_p$ parameterizes the stable reduction of the metacyclic covers of
type ${\bf (x;a)}$. In this section, we want to describe the structure
of $\bar{H}({\bf a})\otimes_R \FF_p$ without proofs. To make the
proofs rigorous one needs to carefully analyze the deformation theory.
Details will appear elsewhere.

\bigskip\noindent {\bf The \'etale case}. In this case all metacyclic
covers of type ${(\bf x;a)}$ have good reduction, Proposition
\ref{maxdim}. Therefore the complete model $\bar{H}({\bf a})$
parameterizes so called admissible covers. It is well known that the
moduli space of admissible covers is smooth, \cite{diss}. We conclude
that $\pi({\bf a}):H({\bf a})\to\PP^1_\lambda$ has good reduction in
this case, confirming Theorem \ref{grthm}.

\bigskip\noindent {\bf The mixed case}. Let $\lambda\in\PP^1_K$ be
such that its reduction $\bar{\lambda}\in\PP^1_k$ is ${\bf
  a}$-ordinary and different from $0,1,\infty$. (See Definition
\ref{ssdef} for the definition of ${\bf a}$-ordinary and ${\bf
a  }$-supersingular.) Let $x_1=0, x_2=1, x_3=\infty, x_4=\lambda$. It
follows from Proposition \ref{genord} that there are $(p-1)/m$
metacyclic covers of type $({\bf x;a})$ with good reduction.  The
other $p(p-1)/m$ such covers have multiplicative bad reduction. If
$\bar{\lambda}$ is ${\bf a}$-supersingular, then all covers have
additive bad reduction.

Suppose $f:Y\to X$ is a metacyclic cover of type $({\bf x;a})$ with
bad reduction. Lemma \ref{vc} implies that $\Xb$ has 5 irreducible
components: the original component $\Xb_0$ and the four primitive
tails $X_1, X_2, X_3, X_4$.  Lemma \ref{vc}.(b) implies that the
ramification invariant $\sigma_i$ corresponding to the tail $X_i$ is
smaller than 1. Therefore it follows from \cite[Lemma
2.3.3]{Rachelfam} that there are only finitely many possibilities for
$\fb|_{X_i}$. This implies that $\bar{\pi}({\bf a}):\bar{H}({\bf
  a})\to\PP^1_R$ is finite, outside the cusps.

It follows from the description of the monodromy given in Section
\ref{grsec} that $\bar{H}({\bf a})\otimes_R \FF_p$ has two irreducible
components $H^\good$ and $H^\bad$. Over the ${\bf a}$-ordinary locus,
$H^\good$ parameterizes metacyclic covers and $H^\bad$ parameterizes
the stable reduction of the covers with multiplicative bad reduction.
The two components intersect above $\lambda\in\Lambda({\bf a})$ (the
${\bf a}$-supersingular $\lambda$'s). 

For $m=2$, the Hurwitz space $H({\bf a})$ is a version of the modular
curve $X_1(p)$ with ordered branch points. (In other words, one adds a
full level 2-structure, cf.\ \cite[Section 2.2]{RRR}.) In this case,
the reduction of $H$ is of course well known.

\bigskip\noindent {\bf The multiplicative case}. All metacyclic covers
of type ${(\bf x;a)}$ have multiplicative bad reduction, Proposition
\ref{mindim}. To describe the reduction of $\pi({\bf a}):H({\bf
  a})\to\PP^1_\lambda$, we first need some results on the stable
reduction of metacyclic covers.

Let ${\bf a}=(a_1, a_2, a_3, a_4)$ be a type such that
$a_1+a_2+a_3+a_4=m$. Let $f:Y\to X$ be a metacyclic cover of type
${\bf a}$, branched at $x_1=0,x_2=1,x_3=\infty, x_4=\lambda$ and
suppose that $\lambda\not\equiv 0,1,\infty\bmod{p}$. Write
$(\Zb_0,\omega)$ for the deformation datum of $f$, as defined in
Definition \ref{defodatum}. As in Remark \ref{typehurw}, one
checks that for $z_i\in\Zb_0$ above $x_i$ we have
\[\gcd(a_i,m)\Ord_{z_i}(\omega)\equiv {a_i}
\pmod{m}.\]
  
  Let $e\in\EE$ be an edge with source $v_0$. Let $x_e$ be the
  corresponding point of $\Xb_0$. Choose a point $z_e\in\Zb_0$ above
  $x_e$ and write $m_e$ for the ramification index of this point in
  $\Zb_0\to\Xb_0$.  Analogous to the notation introduced in Section
  \ref{stablesec}, we define $\nu_e=[h_e/m_e]$ and
  $a_e=m(h_e/m_e-\nu_e)$.
  
  Suppose that the subtree of $T$ with root $e$ contains the primitive
  tail $X_i$. As in Proposition \ref{typehurw}, one
checks that $a_e=a_i$. The analog of Remark \ref{typehurw}
  becomes in this case $\sum_{s(e)=v_0}\nu_e=1$.
It follows that there are two cases.

\begin{itemize}
\item The first possibility is that $\nu_i=0$, for all $i\in\Bp$. In
  this case, there is one new tail $X_5$ which intersects the original
  component.
\item The second possibility is that $\nu_i=1$ for a unique $i\in\Bp$.
\end{itemize}

\begin{defn}\label{specialdef} For $\lambda\in k-\{0,1\}$, let 
  $\Zb_{0}\to \PP^1_k$ be the $m$-cyclic cover of type ${{\bf (x;a)}}$
  branched at $x_1=0, x_2=1, x_3=\infty, x_4=\lambda$. We say that
  $\lambda\in k$ is {\sl $\bf a$-special} if there exists an
  $\omega\in H^0(\Zb_{0}, \Omega^1)_\chi$ and $1\leq i\leq 4$ such
  that $\nu_i(\omega)=1$.  If $\lambda$ is not ${\bf a}$-special, we
  call it is ${\bf a}$-{\sl general}.
\end{defn}

This terminology is inspired by \cite{special}, but the notation of a
special cover in that paper is more restrictive than the way we use
this concept here. 

\begin{lem}\label{numspecial}
  Suppose that $a_1+\cdots +a_4=m$.  There are at most finitely many
  ${\bf a}$-special $\lambda$'s.
\end{lem}

\proof Similar  to \cite[Section 3.5]{special} \Endproof
  
Suppose that $\bar{\lambda}\in\PP^1_k-\{0,1,\infty\}$ is ${\bf
  a}$-general. One checks that $\Xb$ has 6 irreducible components: the
original component $\Xb_0$, the primitive tails $X_1, X_2, X_3, X_4$
and the new tail $X_5$.

\begin{figure}[htb]
\begin{center}
\setlength{\unitlength}{0.0004in}
\begingroup\makeatletter\ifx\SetFigFont\undefined%
\gdef\SetFigFont#1#2#3#4#5{%
  \reset@font\fontsize{#1}{#2pt}%
  \fontfamily{#3}\fontseries{#4}\fontshape{#5}%
  \selectfont}%
\fi\endgroup%
{\renewcommand{\dashlinestretch}{30}
\begin{picture}(4010,2400)(0,-10)
\drawline(1212,2400)(1212,300)
\drawline(2112,2400)(2112,300)
\drawline(3012,2400)(3012,300)
\drawline(12,600)(4212,600)
\drawline(312,2400)(312,300)
\drawline(3912,2400)(3912,300)
\put(4587,600){$\Xb_0$}
\put(162,0){$X_1$ }
\put(1062,0){$X_2$}
\put(1962,0){$X_3$}
\put(2862,0){$X_4$}
\put(3762,0){$X_5$}
\put(-100,1400){$x_1$}
\put(800,1400){$x_2$}
\put(1700,1400){$x_3$}
\put(2600,1400){$x_4$}
\put(312, 1400){\circle*{90}}
\put(1212, 1400){\circle*{90}}
\put(2112, 1400){\circle*{90}}
\put(3012, 1400){\circle*{90}}
\end{picture}
}
\end{center}
\caption{\label{multfig1}Reduction in the multiplicative case for
general $\lambda$}
\end{figure}
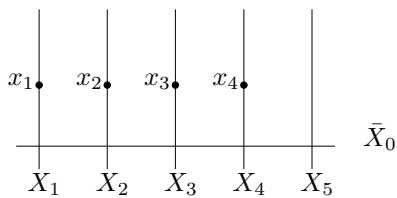

We have $\sigma_i<1$ if $i\leq 4$ and $\sigma_5=2$. We claim that this
implies that, for given $\bar{\lambda}$, that there are only finitely
many possibilities for $\fb$. It follows from \cite[Lemma
2.3.3]{Rachelfam} that there are only finitely many possibilities for
$\fb|_{X_i}$. The intersection point $\tau$ of $X_5$ with $\Xb_0$ is
the image of a zero of $\omega$. Therefore it is obvious that there
are at most finitely many possibilities for $\tau$.

\begin{figure}[htb]
\begin{center}
\setlength{\unitlength}{0.0004in}
\begingroup\makeatletter\ifx\SetFigFont\undefined%
\gdef\SetFigFont#1#2#3#4#5{%
  \reset@font\fontsize{#1}{#2pt}%
  \fontfamily{#3}\fontseries{#4}\fontshape{#5}%
  \selectfont}%
\fi\endgroup%
{\renewcommand{\dashlinestretch}{30}
\begin{picture}(6010,2400)(0,-10)
\drawline(1212,2400)(1212,300)
\drawline(2112,2400)(2112,300)
\drawline(3012,2400)(3012,300)
\drawline(12,600)(3312,600)
\drawline(312,2400)(312,300)
\drawline(-900,1000)(612,1000)
\drawline(-900,1400)(612,1400)
\put(3612,600){$\Xb_0$}
\put(-1400,1400){$X_1$}
\put(-1400,1000){$X_5$} 
\put(162,0){$W$}
\put(1062,0){$X_2$}
\put(1962,0){$X_3$}
\put(2862,0){$X_4$}
\put(-100, 1500){$x_1$}
\put(800,1400){$x_2$}
\put(1700,1400){$x_3$}
\put(2600,1400){$x_4$}
\put(12, 1400){\circle*{90}}
\put(1212, 1400){\circle*{90}}
\put(2112, 1400){\circle*{90}}
\put(3012, 1400){\circle*{90}}
\end{picture}
}
\end{center}
\caption{\label{multfig3}Reduction in the multiplicative  case for ${\bf a}$-special $\lambda$ (Possibility b)}
\end{figure}
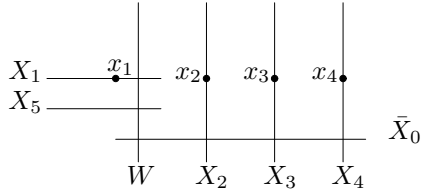

Suppose now that $\bar{\lambda}\in\PP^1_k$ is ${\bf a}$-special. For
simplicity, we suppose that $\nu_1=1$. There are two possibilities for
the stable reduction $\bar{f}$.
\begin{itemize}
\item[(a)] 
 The curve $\Xb$ consists of 5 irreducible components: the
  original component $\Xb_0$ and the primitive tails $X_1, X_2, X_3,
  X_4$, see Figure \ref{multfig2}. 
 The ramification invariants are: $\sigma_1=(m+a_1)/m$ and $\sigma_i=a_i/m$ for $i>1$.
\item[(b)] There is one new tail which does not intersect the original
  component, see Figure \ref{multfig3}.
\end{itemize}

\begin{figure}[htb]
\begin{center}
\setlength{\unitlength}{0.0004in}
\begingroup\makeatletter\ifx\SetFigFont\undefined%
\gdef\SetFigFont#1#2#3#4#5{%
  \reset@font\fontsize{#1}{#2pt}%
  \fontfamily{#3}\fontseries{#4}\fontshape{#5}%
  \selectfont}%
\fi\endgroup%
{\renewcommand{\dashlinestretch}{30}
\begin{picture}(4010,2400)(0,-10)
\drawline(1212,2400)(1212,300)
\drawline(2112,2400)(2112,300)
\drawline(3012,2400)(3012,300)
\drawline(12,600)(3312,600)
\drawline(312,2400)(312,300)
\put(3687,600){$\Xb_0$}
\put(162,0){$X_1$}
\put(1062,0){$X_2$}
\put(1962,0){$X_3$}
\put(2862,0){$X_4$}
\put(-100,1400){$x_1$}
\put(800,1400){$x_2$}
\put(1700,1400){$x_3$}
\put(2600,1400){$x_4$}
\put(312, 1400){\circle*{90}}
\put(1212, 1400){\circle*{90}}
\put(2112, 1400){\circle*{90}}
\put(3012, 1400){\circle*{90}}
\end{picture}
} 
\end{center}
\caption{\label{multfig2}Reduction in the multiplicative  case for ${\bf a}$-special $\lambda$ (Possibility a)}
\end{figure}
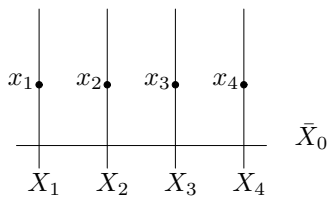 

\Endproof

One can show using a patching argument as in \cite{Henrio} that both
possibilities occur. It follows from \cite[Proposition
2.2.6]{Rachelfam} that there is a one dimensional family of covers of
type (a), which all occur as the reduction of a metacyclic cover.
Therefore the map $\bar{\pi}:\bar{H}\to\PP^1_R$ is not finite: there
are {\sl vertical components} whose image under $\bar{\pi}$ is a
point.  The occurrence of vertical components in the stable reduction
is not unexpected. The same phenomenon occurs for the reduction of the
Hurwitz space of 2-cyclic covers branched at four points to
characteristic 2, \cite[Section 4]{AbrOort98}.

Theorem \ref{grthm} implies that there exists a smooth model of
$H$ over $R$, but the model $\bar{H}$ which has a modular
interpretation is not smooth. 

\section{Examples}\label{exasec}
Suppose that $p> 5$ and choose $2 \neq d|(p-1)/2$. Let
$f:Y\to\PP^1_\CC$ be a cover with Galois group $\PSL_2(p)$ branched at
three points of order $p,p,d$.  In this section we show that there
exists an $m|(p-1)$ and a type ${\bf a}$ such that $f$ is isomorphic
to the quotient of $\varpi({\bf a}):\HH({\bf a})\to\PP^1$ by the
center of its Galois group. A similar result is shown in \cite[Section
2.2]{RRR} for $\PSL_2(p)$-covers of $\PP^1_\CC$ branched at three
points of order $p$. In that case one may take $m=2$ and ${\bf
  a}=(1,1,1,1)$.  The cover $\varpi({\bf a})$ is in this case
isomorphic to $X(2p)\to X(2)$, where $X(N)$ is the modular curve with
full level-$N$ structure.

Recall that $\SL_2(p)$ has two conjugacy classes of elements of order
$p$, which we denote by $pA$ and $pB$, \cite[Lemma 3.27]{Voelklein95}.
We make the convention that
\[\left(\begin{array}{cc}1&1\\0&1\end{array}\right)\in pA.\]
Choose a primitive $(p-1)$th root of unity $\xi\in\FF_p^\times$.  For
$0<i<(p-1)/2$, let $\C(i)=\{A\in\SL_2(p)|\, \Tr(A)=\xi^i+\xi^{-i} \}.$
Then $\C(i)$ is a conjugacy class of $\SL_2(p)$ of elements whose
order is equal to the order of $\xi^i$ in $\FF_p^\times$. Moreover,
every element of $\SL_2(p)$ of order $n$ dividing $p-1$ with $n>2$ is
contained in some $\C(i)$.   For a triple of conjugacy classes ${\bf
  C}=(C_1, C_2, C_3)$ of $\SL_2(p)$, we write $\Niin_3({\bf C})$ for
the set of isomorphism classes of $\SL_2(p)$-cover $Y\to\PP^1_\CC$
with class vector ${\bf C}$.  (This means that the canonical generator
of some point of $Y$ above $x_i$ is contained in the conjugacy class
$C_i$, with respect to some fixed system of roots of unity.)

\begin{lem}
  Let $0<i<(p-1)/2$. Then either $|\Niin_3(pA,pA,\C(i))|=0$ and
  $|\Niin_3(pA,pB,\C(i))|=1$, or $|\Niin_3(pA,pA,\C(i))|=1$ and
  $|\Niin_3(pA,pB,\C(i))|=0$.
\end{lem}

\proof
The proof is similar to \cite[Lemma 3.29]{Voelklein}.
\Endproof

The lemma says that the triples $(pA,pA,\C(i))$ and $(pA,pB, \C(i))$
are {\sl rigid}.  Since there is an outer automorphism of $\SL_2(p)$
which interchanges the conjugacy classes $pA$ and $pB$ and fixes
$\C(i)$, the lemma may be rephrased as follows. For every
$0<i<(p-1)/2$, there exists a unique $\SL_2(p)$-cover of $\PP^1_\CC$
with class vector $(p\ast,p\ast\ast, \C(i))$ {\sl up to isomorphism}.

We will now construct types ${\bf a}$ such that $\varpi({\bf a})$ is
ramified of order $p,p,d$, for some $d| (p-1)$. The cover $\varpi({\bf
  a})$ is branched at $0$ and 1 of order $p$ if and only if
$a_1+a_4=a_2+a_4=m$, by Proposition \ref{rampi}. Therefore the only
types we need to consider are types of the form $(a,a,m-a, m-a)$, for
some $m|(p-1)$ with $\gcd(a,m)=1$.

\begin{lem} Suppose  $m\neq 2$ divides $p-1$. Choose $0<a<m$ such that
  $\gcd(m,a)=1$ and write ${\bf a}=(a,a,m-a,m-a)$. Then the Galois
  group of $\varpi({\bf a})$ is $\SL_2(p)$ if $m$ is odd and
  $\PSL_2(p)$ if $m$ is even.
\end{lem} 

\proof We first note that $p$ divides the order of the Galois group of
$\varpi({\bf a})$, since the ramification index at $0$ and $1$ is $p$.
Therefore we are not in the exceptional case, Definition
\ref{exceptionaldef}. 

We denote by $B_1, B_2, B_3$ the matrices defined in
(\ref{B1})--(\ref{B3}). Then $\det(B_1)=\det(B_2)=1$ and
$\det(B_3)=(\xi^{\alpha a})^2$. Corollary \ref{Galpi} implies that the
Galois group of $\varpi({\bf a})$ is $\SL_2(p)$ if $m$ is odd and
$\PSL_2(p)$ if $m$ is even. \Endproof

For every type ${\bf a}$, we can define an $\SL_2(p)$-cover
$\varpi'({\bf a})$. If $m$ is odd, this is just $\varpi({\bf a})$. If
$m$ is even, we choose a lift of $\varpi({\bf a})$ to an
$\SL_2(p)$-cover of $\PP^1$ branched at three points. We may suppose
that $\varpi'({\bf a})$ is branched of order $p$ at 0 and 1. This
uniquely characterizes the lift.

The following
proposition states that there exists a type ${\bf a}$ such that
$\varpi'({\bf a})$ has class vector $(C_1,C_2, \C(i))$, for some choice
of $C_1, C_2$. In other words, every $\SL_2(p)$-cover with class
vector $(p*, p*, \C(i))$ is isomorphic to one of the covers
$\varpi'({\bf a})$. 

\begin{prop}\label{classvec}
  Let $2\neq m|(p-1)$ and ${\bf a}=(a,a,m-a, m-a)$. Write
  $\alpha=(p-1)/m$.  Suppose $a<m/2$ and $\gcd(m,a)=1$. Then the class
  vector of $\varpi'({\bf a}):\HH({\bf a})\to\PP^1$ is $(C_1, C_2,
  \C(a\alpha))$, for some $C_1, C_2\in\{pA, pB\}$.
\end{prop}
  
\proof The Nielsen tuple of $\varpi'({\bf a})$ is $(B_1, B_2, B_3')$
with $B'_3=B_3/\xi^{\alpha a}\in\SL_2(p)$.  Therefore
$\Tr(B'_3)=\xi^{\alpha a}+\xi^{-\alpha a}$ and $B'_3\in\C(\alpha a)$.
\Endproof

\begin{cor}\label{isocor}
  Suppose $d\neq 2$ divides $(p-1)/2$.  Let
  $f:Y\to\PP^1_\CC$ be a cover with Galois group $\PSL_2(p)$ which is
  branched only at $0,1,\infty$ with ramification index $p,p,d$,
  respectively. Then there exists a type ${\bf a}$ such that $f$ is
  isomorphic to the quotient of $\varpi({\bf a})\otimes_K \CC$ by the
  center of the Galois group $G$ of $\varpi({\bf a})$.
\end{cor}

\proof This follows from Proposition \ref{classvec}.
\Endproof

One may expect Corollary \ref{isocor} to hold more generally. Let
$G=\SL_2(p)$ and suppose that ${\bf C}=(C_1,C_2,C_3)$ is a triple of
conjugacy classes of $G$.  In general, ${\bf C}$ is not rigid, but it
is still {\sl linearly rigid}, \cite[Example 2.4]{StVo}. This means
the following. Let $(M_1, M_2, M_3)$ be a triple of matrices in
$\GL_2(p)$ such that $M_i\in C_i$ and $M_1M_2M_3=1$. Then if $(M_1',
M_2', M_3')$ is another such triple, there exists an $N\in\GL_2(p)$
such that $N M_i' N^{-1}=M_i$ for all $i$. This was essentially
already known to Riemann, see \cite[Introduction]{Katzrigid}.

It follows from Proposition \ref{rampi} that the ramification indices
of $\varpi({\bf a})$ are either $p$ or divide $p-1$. Therefore there
exist covers of $\PP^1$ branched at three points with Galois group
$\PSL_2(p)$, $\SL_2(p)$ or $\PGL_2(p)$ which are not isomorphic to
some $\varpi({\bf a})$. To obtain covers with ramification indices
dividing $p+1$ one should omit the condition $p\equiv 1\bmod{m}$.

\bigskip\noindent {\bf Comparison to Raynaud's criterion for good
  reduction}.
In the rest of this paper, we want to compare Theorem \ref{grthm} to
Raynaud's criterion for good reduction, \cite{Raynaud98}. We
illustrate that our result does not follow from Raynaud's Criterion.

Suppose that $m=5$ and $p\equiv 1\bmod{5}$. Choose a primitive 5th
root of unity $\zeta_5\in\FF_p^\times$. For $1\leq i\leq 4$, we write
${\bf a}_i=(i,i,i,-3i)$. One checks that ${\bf a}_i$ is
non-exceptional. Therefore $\varpi_i:=\varpi({\bf a}_i)$ is a Galois cover
with Galois group $\SL_2(p)$, Corollary \ref{Galpi}.  Proposition
\ref{rampi} implies that $\varpi_i$ is branched at three points of
order 5.

Let $5A$ (resp.\ $5B$) be the conjugacy class of $\SL_2(p)$ of
matrices with trace $\zeta_5+ \zeta_5^4$ (resp.\ $\zeta_5^2+\zeta_5^3$).
Write ${\bf C}_A=(5A, 5A, 5A)$ and ${\bf C}_B=(5B, 5B, 5B)$. As
before, we denote by $\Niin_3({\bf C})$ the  Nielsen class.

\begin{lem}\label{555} 
\begin{itemize}
\item[(a)] We have $\varpi_1, \varpi_4\in \Niin_3({\bf C}_B)$ and 
 $\varpi_2, \varpi_3\in \Niin_3({\bf C}_A)$.
\item[(b)] The sets $\Niin_3({\bf C}_A)$ and $\Niin_3({\bf C}_B)$
  both contain  exactly two elements.
\end{itemize}
\end{lem}

\proof For each type, we defined matrices $B_1, B_2, B_3$ in
(\ref{B1})--(\ref{B3}). We denote these matrices corresponding to the
type ${\bf a}_i$ by $B_1(i), B_2(i), B_3(i)$.

The class vector of $\varpi_i$ is described by matrices
$B_1'(i), B_2'(i), B_3'(i)$ which are characterized by the following
two properties:
\begin{itemize}
\item $B_j'(i)=B_j(i)\cdot \gamma_j(i) I$, for some scalar matrix
      $\gamma_j(i) I\in M$,  such that $\det(B_j'(i))=1$,
\item  the order of $B_j'(i)$ is 5.
\end{itemize}
It follows that $\gamma_j(i)=\zeta_5^{-i}$, for every $j$. Part (a)
follows from this.  The triples ${\bf C}_A$ and ${\bf C}_B$ are
linearly rigid, \cite[Example 2.4]{StVo}. The outer automorphism group
of $\SL_2(p)$ has order two and fixes the conjugacy classes of order
five. Therefore, the Nielsen class has either zero or two elements and
(b) follows from (a). \Endproof

In $\SL_2(p)$ there are two conjugacy classes of elements of order 5.
Therefore the triples ${\bf C}_A$ and ${\bf C}_B$ are rational over
$\QQ(\sqrt{5})$, \cite[Section 7.1]{Serretopics}.  Lemma \ref{555}.(b)
implies that the $\SL_2(p)$-covers $\varpi_i$ are defined over some
extension $K_i$ of $\QQ(\sqrt{5})$ of degree at most two.  Let $\wp$
be a prime of $K_i$ above $p$ and write  $e_i$ for the ramification index of $\wp$. Since the degree
of $K_i$ over $\QQ(\sqrt{5})$ is at most two, we have $e_i\leq 2$.
Raynaud's Criterion for good reduction states in this case that if
$\varpi_i$ has bad reduction then $e_i$ is greater than or equal to the
number of Sylow $p$-subgroups in $\SL_2(p)$, i.e.\ $e_i\geq 2$.
Therefore we cannot conclude directly whether $\varpi_i$ has good
reduction or not, unless we have more information on the field of
definition $K_i$ of $\varpi_i$. But such information is hard to obtain
for triples which are not rigid. In what follows we take an
alternative approach. We use Theorem \ref{grthm} to decide whether
$\varpi_i$ has good reduction and deduce from this information on the
field of definition $K_i$.

 The following proposition is essentially the
 statement of Raynaud's criterion for good reduction, in our special
 case.

\begin{prop}\label{Racrit}
  The cover $\varpi_i$ has good reduction at $\wp$ if and only if
  $e_i=1$.
\end{prop}

\proof If $\varpi_i$ has good reduction at $\wp$ it follows from
Beckmann's Theorem that $e_i=1$, \cite{Beckmann89}.  Suppose that
$\varpi_i$ has bad reduction. Raynaud's Criterion \cite{Raynaud98}
implies that $e_i\geq 2$. In \cite{Raynaud98} the assumption is made
that the Galois group $G$ has trivial center, but this is not
essential.  \Endproof

\begin{prop}
\begin{itemize}
\item[(a)] Let $f:Y\to\PP^1$ be an $\SL_2(p)$-cover with class vector
  ${\bf C}_A$. Then $f$ has bad reduction.
\item[(b)] Let $f:Y\to\PP^1$ be an $\SL_2(p)$-cover with class vector
  ${\bf C}_B$. Then $f$ has good reduction.
\end{itemize}
\end{prop}

\proof Lemma \ref{555} implies that, over $\CC$, there are two covers
with class vector ${\bf C}_A$ and two with class vector ${\bf C}_B$.
Write $\Niin_4({\bf C}_B)=\{f_1, f_4\}$ and $\Niin_4({\bf C}_A)=\{f_2,
f_3\}$.  Let $K/\QQ$ be a field over which all $f_i$ are defined. The
definition of the type ${\bf a}_i$ depends on a choice of primitive
$5$th root of unity $\xi_5\in K$.  Recall that we have also fixed a
primitive $5$th root of unity $\zeta_5\in\FF_p^\times$.  Choose a
prime ideal $\wp$ of $K$ above $p$ such that $\xi_5\equiv \zeta_5
\bmod{\wp}$.

We may assume that $f_1=\varpi_1$ and $f_2=\varpi_2$.  Theorem
\ref{grthm} implies that $f_1$ has good reduction and $f_2$ has bad
reduction at $\wp$.  It follows from Proposition \ref{Racrit} that
$f_2$ is defined over a field $K_2$ such that $\QQ(\sqrt{5})\subset
K_2\subset K$ in which the ramification index of $\wp\cap K_2$ is two.
Let $\sigma_2$ be a generator of $\Gal(K_2, \QQ(\sqrt{5}))$. Then
$f_3=f_2^\sigma$. we conclude that $f_3$ has bad reduction at $\wp$
also.

The Galois group of $\QQ(\sqrt{5})/\QQ$ permutes
  the class vectors ${\bf C}_A$ and ${\bf C}_B$, therefore the covers
  $f_1,\, f_2,\, f_3,\,f_4$ are conjugated over $K$ by an automorphism
  of order four. It follows that we may take $K$ to be the Galois closure
  of $K_2/\QQ$. The Galois group of $K/\QQ$ is a transitive group on
  four letters which contains an element of order four and has a
  quotient of order two, i.e.\ $\Gal(K,\QQ)$ is either cyclic of order
  four or a dihedral group of order eight. 
  
  Suppose $\Gal(K,\QQ)\simeq \ZZ/4$. Then $K_2/\QQ$ is Galois and
  $f_1$ is defined over $K_2$. Moreover, $f_1^\sigma=f_4$. This
  implies that $K_2$ is a minimal field of definition for $f_1$. Since
  $\wp$ is ramified in $K_2$, Proposition \ref{Racrit} implies that
  $f_1$ has bad reduction at $\wp$ contradicting Theorem \ref{grthm}.
  We conclude that $\Gal(K,\QQ)$ is a dihedral group of order eight.
  We represent the elements of $\Gal(K,\QQ)$ as permutations on
  $\{1,2,3,4\}$. It follows that $K_2$ is the subfield of $K$ of
  invariants under $(1\,4)$. 
  
  The covers $f_1$ and $f_4$ may be defined over the subfield $K_1$ of
  $K$ of invariants of $(2\,3)$. The permutation $(1\,4)(2\,3)$ restricts
  to an automorphism of $K_1$ which permutes $f_1$ and $f_4$,
  therefore $K_1$ is a minimal field of definition of $f_1$. Since
  $f_1$ has good reduction at $\wp$, we conclude that $\wp\cap K_1$ is
  unramified in $K_1$. Therefore $f_4$ has good reduction at $\wp$
  also. 
  
  We remark that the choice of $\wp$ depends on the choice of the 5th
  root of unity $\xi_5\in K$. The permutation
  $(1\,2\,4\,3)\in\Gal(K,\QQ)$ sends $\wp$ to a different prime $\wp'$
  and also changes the class vector of a cover $f_i$.  If $f_i$ has
  good reduction at $\wp$ then $f_i$ has bad reduction at $\wp'$, and
  conversely.  \Endproof

\bibliographystyle{abbrv} \bibliography{hurwitz} 

\begin{thebibliography}{10}

\bibitem{AbrOort98}
D.~Abramovich and F.~Oort.
\newblock Stable maps and {H}urwitz schemes in mixed characteristic.
\newblock In {\em Advances in algebraic geometry motivated by physics}, number
  276 in Contemp.\ Math., pages 89--100. Amer.\ Math.\ Soc., 2001.

\bibitem{Beckmann89}
S.~Beckmann.
\newblock Ramified primes in the field of moduli of branched coverings of
  curves.
\newblock {\em J.\ of Algebra}, 125:236--255, 1989.

\bibitem{Berger}
G.~Berger.
\newblock Fake congruence subgroups and the {H}urwitz monodromy group.
\newblock {\em J.\ Math.\ Sci.\ Univ.\ Tokyo}, 6:559--574, 1999.

\bibitem{thesis}
I.~I. Bouw.
\newblock The $p$-rank of ramified covers of curves.
\newblock {\em Compositio Math.}, 126:295--322, 2001.

\bibitem{RRR}
I.~I. Bouw and R.~J. Pries.
\newblock Rigidity, reduction, and ramification.
\newblock To appear in {\it Math. Ann.}

\bibitem{modular}
I.~I. Bouw and S.~Wewers.
\newblock Reduction of covers and {H}urwitz spaces.
\newblock math.AG/0005120.

\bibitem{CW}
C.~Chevalley and A.~Weil.
\newblock \"{U}ber das {V}erhalten der {I}ntegrale erster {G}attung bei
  {A}utomorphismen des {F}unktionenk\"orpers.
\newblock {\em Abh.\ Math.\ Sem.\ Univ.\ Hamburg}, 10:358--361, 1934.

\bibitem{CoMc}
R.~Coleman and W.~McCallum.
\newblock Stable reduction of {F}ermat curves and {J}acobi sum {H}ecke
  characters.
\newblock {\em J.\ Reine Angew.\ Math.}, 385:41--101, 1988.

\bibitem{DelRap}
P.~Deligne and M.~Rapoport.
\newblock Les sch\'emas de modules de courbes elliptiques.
\newblock In {\em Modular functions of one variable II}, number 349 in LNM,
  pages 143--316. Springer-Verlag, 1972.

\bibitem{Henrio}
Y.~Henrio.
\newblock Arbres de {H}urwitz et automorphismes d'ordre $p$ des disques et
  couronnes $p$-adiques formels.
\newblock Phd-thesis, University of Bordeaux, 1999.

\bibitem{Huppert}
B.~Huppert.
\newblock {\em Endliche Gruppen I}.
\newblock Number 134 in Grundlehren. Springer-Verlag, 1967.

\bibitem{Ihara74}
Y.~Ihara.
\newblock On the differentials associated to congruence relations and the
  {S}chwarzian equations defining uniformizations.
\newblock {\em J.\ Fac.\ Sci.\ Univ.\ Tokyo Sect.\ IA Math.}, 21:309--332,
  1974.

\bibitem{Katzrigid}
N.~Katz.
\newblock {\em Rigid local systems}.
\newblock Number 139 in Annals of Mathematics Studies. Princeton Univ.\ Press,
  1996.

\bibitem{KatzMazur}
N.~M. Katz and B.~Mazur.
\newblock {\em Arithmetic moduli of elliptic curves}.
\newblock Number 108 in Annals of Mathematics Studies. Princeton Univ.\ Press,
  1985.

\bibitem{Claus}
C.~Lehr.
\newblock Reduction of $p$-cyclic covers of the projective line.
\newblock {\em Manuscripta Math.}, 106:151--175, 2001.

\bibitem{Matignon01}
M.~Matignon.
\newblock Vers un algorithme pour la r\'eduction stable des rev\^etements
  $p$-adiques de la droite projective sur un corps $p$-adique.
\newblock math.NT/0112042.

\bibitem{Milne}
J.~S. Milne.
\newblock {\em \'{E}tale cohomology}.
\newblock Princeton Univ.\ Press, 1980.

\bibitem{Rachelfam}
R.~J. Pries.
\newblock Families of wildly ramified covers of curves.
\newblock To appear in {\it Amer.\ J.\ of Math.}

\bibitem{Raynaudfest}
M.~Raynaud.
\newblock $p$-groupes et r\'eduction semi-stable des courbes.
\newblock In P.~Cartier, editor, {\em Grothendieck festschrift III}, number~88
  in Progress in Math., pages 179--197. Birkh\"auser, 1990.

\bibitem{Raynaud94}
M.~Raynaud.
\newblock Rev\^etement de la droite affine en caract\'eristique $p>0$ et
  conjecture d'{A}bhyankar.
\newblock {\em Invent.\ Math.}, 116:425--462, 1994.

\bibitem{Raynaud98}
M.~Raynaud.
\newblock Sp\'ecialisation des rev\^etements en caract\'eristique $p>0$.
\newblock {\em Ann.\ Sci.\ \'Ecole Norm.\ Sup.}, 32:87--126, 1999.

\bibitem{Saidi00}
M.~Sa{\"\i}di.
\newblock Wild ramification and a vanishing cycles formula.
\newblock math.AG/0106248.

\bibitem{SerreMex}
J.-P. Serre.
\newblock Sur la topologie des vari\'et\'es alg\'ebriques en caract\'eristique
  $p$.
\newblock {\em Symb.\ Int.\ Top.\ Alg.\ Mexico}, pages 24--53, 1958.
\newblock (\OE vres no.\ 38).

\bibitem{Serretopics}
J.-P. Serre.
\newblock {\em Topics in Galois theory}.
\newblock Number~1 in Research notes in mathematics. Jones and Bartlett
  Publishers, 1992.

\bibitem{StVo}
K.~Strambach and H.~V\"olklein.
\newblock On linearly rigid triples.
\newblock {\em J.\ reine angew.\ Math.}, 510:57--62, 1999.

\bibitem{Voelklein93}
H.~V{\"o}lklein.
\newblock Braid group action via $\gl_n(q)$ and $\uu_n(q)$, and {G}alois
  realizations.
\newblock {\em Israel J.\ Math.}, 82:405--427, 1993.

\bibitem{Voelklein95}
H.~V{\"o}lklein.
\newblock Cyclic covers of $\pp^1$ and {G}alois action on their division
  points.
\newblock In M.~D. Fried, editor, {\em Recent developments in the inverse
  {G}alois problem}, number 186 in Contemp.\ Math., pages 91--107, 1995.

\bibitem{Voelklein}
H.~V{\"o}lklein.
\newblock {\em Groups as {G}alois groups}.
\newblock Number~53 in Cambridge Studies in Adv.\ Math. Cambridge Univ. Press,
  1996.

\bibitem{special}
S.~Wewers.
\newblock Reduction and lifting of special metacyclic covers.
\newblock math.AG/0105052, to appear in {\it Ann.\ Sci.\ \'Ecole Norm.\ Sup.}

\bibitem{Stefan01}
S.~Wewers.
\newblock Three point covers with bad reduction.
\newblock In preparation.

\bibitem{diss}
S.~Wewers.
\newblock Construction of {H}urwitz spaces.
\newblock Thesis, Preprint No.\ 21 of the IEM, Essen, 1998.

\bibitem{Yoshida}
M.~Yoshida.
\newblock {\em Hypergeometric functions, my love}.
\newblock Number E 32 in Aspects of Math. Vieweg, 1997.

\end{thebibliography}
\end{document}